\documentclass[11pt]{article}
\usepackage{amsmath,amssymb,amsthm}
\usepackage[a4paper,margin=2.5cm]{geometry}
\usepackage{graphicx}
\usepackage{float}
\usepackage[hidelinks]{hyperref}

\numberwithin{equation}{section}
\newtheorem{lemma}{Lemma}[section]
\newtheorem{proposition}[lemma]{Proposition}
\newtheorem{theorem}[lemma]{Theorem}
\newtheorem{corollary}[lemma]{Corollary}
\theoremstyle{definition}
\newtheorem{definition}[lemma]{Definition}
\theoremstyle{remark}
\newtheorem{remark}[lemma]{Remark}

\newcommand{\R}{\mathbb{R}}

\newcommand{\E}{\mathbb{E}}
\newcommand{\PP}{\mathbb{P}}
\newcommand{\cH}{\mathcal{H}}
\newcommand{\D}{\mathbb{D}}

\title{The stochastic phase field crystal equation with conservative space-time white noise: well-posedness, regularity and existence of a density}
\author{Atef Lechiheb}
\date{}

\begin{document}
\maketitle
{\renewcommand{\thefootnote}{}\footnotetext{\noindent Atef Lechiheb, Toulouse School of Economics, Universit\'e Toulouse Capitole, Toulouse, France.\\
E-mail: \texttt{atef.lechiheb@tse-fr.eu}\quad ORCID: 0000-0002-9896-3408}}
\begin{abstract}
We study the phase field crystal equation $\partial_tu=\Delta[(1+\Delta)^2u+u^3-\varepsilon u]+\sigma\partial_x\dot W$ on $[0,\pi]$ with the boundary conditions $\partial_xu=\partial_x^3u=\partial_x^5u=0$ at $x=0,\pi$, perturbed by a conservative space-time white noise. This sixth order equation is the standard model for crystal growth at atomic length scales and diffusive time scales, and the conservative noise is the form of thermal fluctuations used in the physics literature. We prove that for every initial condition in $H^1$ there exists a unique global mild solution, continuous in space and time, which preserves the mean and whose supremum has finite moments of all orders. We determine its H\"older regularity: away from $t=0$, it is H\"older continuous in time with any exponent smaller than $1/4$, and $C^1$ in space with a spatial derivative which is H\"older continuous with any exponent smaller than $1/2$. Finally, using Malliavin calculus and a localization argument, we prove that for $\sigma\neq0$ the law of $u(t,x)$ is absolutely continuous with respect to the Lebesgue measure for every $t>0$ and every $x$. A specific feature of the model is a second neutral mode of the linear operator, besides the constant mode, which requires a separate treatment in the a priori estimates.
\end{abstract}



\section{Introduction}

\subsection{The model}

The phase field crystal (PFC) model was introduced by Elder, Katakowski, Haataja and Grant \cite{EKHG02,EG04} to describe the evolution of crystalline materials at atomic length scales and on diffusive time scales. The unknown $u$ is related to the local atomic density. In the solid phase it is periodic, with the symmetry of the crystal lattice, while it is constant in the liquid phase. The model is the conserved gradient flow, in $H^{-1}$, of the Swift-Hohenberg type free energy
\[
\mathcal E(u)=\int_D\Big(\frac12\big((1+\Delta)u\big)^2+\frac14u^4-\frac\varepsilon2u^2\Big)dx,
\]
which leads to the sixth order equation
\[
\partial_tu=\Delta\big[(1+\Delta)^2u+u^3-\varepsilon u\big].
\]
Compared with phase field models of Cahn-Hilliard type, this approach couples elasticity and plasticity to the thermodynamics without auxiliary fields, and it has been used for grain growth, dendritic and eutectic solidification, epitaxial growth and defect dynamics.

Thermal fluctuations are usually added through a conservative noise, that is, a Gaussian forcing $\eta$ with covariance $\langle\eta(x,t)\eta(x',t')\rangle=-\alpha\Delta\delta(x-x')\delta(t-t')$ \cite{EKHG02,XMYL20}. Such a term is needed, for instance, to obtain nucleation from a homogeneous liquid state. Whether a noise term should be included at all has been discussed in the derivation of the model from dynamical density functional theory: Archer and Rauscher \cite{AR04} showed that the noise is justified when the density is understood as coarse grained in time. In one space dimension, this conservative noise is $\sigma\partial_x\dot W$, where $\dot W$ is a space-time white noise, and we are led to
\begin{equation}\label{eq:intro}
\partial_tu=\Delta\big[(1+\Delta)^2u+f(u)\big]+\sigma\partial_x\dot W,\qquad f(u)=u^3-\varepsilon u,
\end{equation}
on $D=[0,\pi]$ with the boundary conditions $\partial_xu=\partial_x^3u=\partial_x^5u=0$ at $x=0$ and $x=\pi$. These three conditions are needed to define a realization of the sixth order operator; we refer to them as Neumann conditions.

\subsection{Related works}

The deterministic PFC equation and its modified (damped wave) version are well understood. Energy stable and convergent numerical schemes were constructed by Wise, Wang and Lowengrub \cite{WWL09}, global smooth solutions of the modified equation were obtained by Wang and Wise \cite{WW10}, and well-posedness and attractors were studied by Grasselli and Wu \cite{GW14}.

On the stochastic side, the closest equation is the stochastic Cahn-Hilliard equation, for which the literature is extensive. Da Prato and Debussche \cite{DPD96} proved existence and uniqueness with additive space-time white noise, and Cardon-Weber \cite{CW01} proved the existence of a density by Malliavin calculus, a method that we follow here. Cardon-Weber and Millet \cite{CWM04} developed a general framework for SPDEs driven by operators of arbitrary order which are parabolic in the sense of Petrovski\u{\i}, with globally Lipschitz coefficients, and combined it with a priori estimates specific to the Cahn-Hilliard equation to treat that case. The case of fractional noise in time was studied by Bo, Jiang and Wang \cite{BJW08} for existence (see \cite{Le26} for a discussion of one of the estimates of that paper) and by Dimitriou, Farazakis and Karali \cite{DFK26} for the density, and Antonopoulou, Farazakis and Karali \cite{AFK18} studied the Malliavin calculus for a combined Cahn-Hilliard/Allen-Cahn equation. For phase field models, Bertini, Brassesco, Butt\`a and Presutti \cite{BBBP02} proved global existence and uniqueness for a one dimensional system of Caginalp type on $\R$, which couples a stochastic Allen-Cahn equation for the order parameter $m$ with a heat equation for the temperature $h$; there the conserved quantity is $m+h$ and the space-time white noise acts on the order parameter only, in a non-conservative way. Our equation is of a different nature: it is a single equation of order six, the noise is conservative, and the linear operator has two neutral modes. The stochastic Swift-Hohenberg equation, which shares the operator $(1+\Delta)^2$ with the PFC equation but has no conservative structure, has been studied mainly near a change of stability, through amplitude and modulation equations \cite{MBK13,BN22}.

For the stochastic PFC equation itself, the works we are aware of are numerical, see for instance \cite{XMYL20} and the references therein. To the best of our knowledge, no well-posedness result is available for \eqref{eq:intro}, and the general results of \cite{CWM04} do not apply directly, since the nonlinearity $\Delta f(u)$ is not globally Lipschitz and the a priori estimates have to be derived for the specific structure of the PFC equation, with its two neutral modes.

\subsection{Main results and difficulties}

Our main results are the following, for $u_0\in H^1(D)$ and $\varepsilon,\sigma\in\R$.
\begin{enumerate}
\item \emph{Global well-posedness} (Theorem~\ref{thm:global}). Almost surely, \eqref{eq:intro} has a unique global mild solution $u\in C([0,T]\times D)$, $\int_Du(t,x)\,dx$ is constant in time, and $\sup_{t,x}|u(t,x)|$ has finite moments of all orders.
\item \emph{H\"older regularity} (Theorem~\ref{thm:regularity}). On $[t_0,T]\times D$, $t_0>0$, the solution is $\gamma_1$-H\"older in time for all $\gamma_1<\frac14$ and Lipschitz in space; moreover $u(t,\cdot)\in C^1(D)$ for $t>0$, and $\partial_xu$ is $\gamma$-H\"older in space for all $\gamma<\frac12$.
\item \emph{Existence of a density} (Theorem~\ref{thm:density}). If $\sigma\neq0$, the law of $u(t,x)$ is absolutely continuous for every $t\in(0,T]$ and $x\in D$.
\end{enumerate}

Three features distinguish the PFC equation from the Cahn-Hilliard case. First, the linear operator $\Delta(1+\Delta)^2$ has eigenvalues $\lambda_k=k^2(1-k^2)^2$ and therefore two neutral modes, $k=0$ and $k=1$. The constant mode is handled by the conservation of the mean, as for Cahn-Hilliard, but the mode $\cos x$ is neither damped by the linear part nor fixed by a conservation law, and it requires a separate argument in the a priori estimate (Lemma~\ref{lem:interp}). Second, the operator is of order six, which gives stronger smoothing in the deterministic convolution, of order $t^{-1/2}$ instead of $t^{-3/4}$ for the term $\Delta_yG$. Third, the conservative noise costs one derivative, while the operator of order six gains two derivatives with respect to Cahn-Hilliard. As a result, the stochastic convolution is $C^1$ in space (Lemma~\ref{lem:Zx}) and $\gamma$-H\"older in time for $\gamma<\frac14$. For comparison, with the same conservative noise the Cahn-Hilliard convolution is only $\gamma$-H\"older in space for $\gamma<\frac12$ and in time for $\gamma<\frac18$; with an additive space-time white noise it has the same spatial regularity as ours but is smoother in time, with exponents $\gamma<\frac38$ \cite{CW01}. This last property is what makes our a priori estimate work: in the $L^2$ estimate for $v=u-Z$, which we derive on a spectral Galerkin approximation in order to make every step rigorous, the cubic part of the nonlinearity has a good sign and absorbs the terms involving $\partial_xZ$. Combined with an $H^{-1}$ estimate which exploits the quartic potential (Lemma~\ref{lem:Hminus1}), this yields a priori bounds that are polynomial in the noise, and hence moments of all orders for $\sup_{t,x}|u(t,x)|$.

\subsection{Applications}

The noise in \eqref{eq:intro} is not a technical addition: several phenomena described by the PFC model cannot occur without it. The deterministic dynamics started from a homogeneous liquid state remains homogeneous, and thermal fluctuations are needed to trigger crystallization \cite{AR04,TTPTG10}. Once crystallites have appeared at different places with different orientations, their growth and their encounters produce the grain boundaries of polycrystalline materials, whose dynamics is a main application of the model \cite{EG04,PE10}. Fluctuations also produce the capillary waves at the crystal-liquid interface and influence the formation and motion of defects \cite{PE10}, and the same equation is used for the periodic patterns of block copolymers \cite{XMYL20}. All these simulations rely on equations of the form \eqref{eq:intro}, for which our results provide a rigorous basis: well-posedness, which is needed to give a meaning to the simulated paths and to study the convergence of numerical schemes; regularity, which determines the rates that such schemes can reach; and the existence of a density, which is the starting point for quantitative statements on the probability that $u(t,x)$ exceeds a given level. We illustrate these points numerically in Section~\ref{sec:numerics}, where the noise is seen to drive the solution away from the homogeneous state towards one of two periodic profiles, chosen at random.

\subsection{Organization of the paper}

Section~2 presents the spectral setting, the Green function and the stochastic convolution. Section~3 proves local existence and uniqueness by a pathwise fixed point argument. Section~4 contains the a priori estimate and the global existence. Section~5 is devoted to the H\"older regularity of the solution, and Section~6 to its Malliavin differentiability and the existence of a density. Section~\ref{sec:numerics} illustrates these results on numerical simulations. Appendix~\ref{app:tools} collects a Gronwall lemma with weakly singular kernel, elementary estimates on spectral sums, and two classical results on Gaussian fields that are used throughout.
\section{Setting, Green function and stochastic convolution}\label{sec:setting}

\subsection{The model and the noise}

Let $D=[0,\pi]$, $T>0$, and let $(\Omega,\mathcal F,(\mathcal F_t)_{t\in[0,T]},\PP)$ be a filtered probability space satisfying the usual conditions, carrying a Brownian sheet $W=\{W(t,x),(t,x)\in[0,T]\times D\}$ adapted to $(\mathcal F_t)$. We write $W(ds,dy)$ for the associated martingale measure and use the stochastic integral of Walsh \cite{Wa86}: for a deterministic $h\in\cH:=L^2([0,T]\times D)$, $\int_0^T\!\int_D h(s,y)\,W(ds,dy)$ is a centred Gaussian variable with variance $\|h\|^2_{\cH}$. Equivalently, the time derivative of $W$ is a space-time white noise, or the derivative of a cylindrical Wiener process on $L^2(D)$; the two points of view coincide on $D$.

We consider the stochastic phase field crystal equation with conservative noise
\begin{equation}\label{eq:pfc}
\partial_t u = \Delta\big[(1+\Delta)^2 u + f(u)\big] + \sigma\, \partial_x \dot W, \qquad (t,x)\in(0,T]\times D,
\end{equation}
with $f(u)=u^3-\varepsilon u$, $\varepsilon\in\R$, $\sigma\in\R$, initial condition $u(0,\cdot)=u_0$, and the Neumann boundary conditions
\begin{equation}\label{eq:bc}
\partial_x u=\partial_x^3 u=\partial_x^5 u=0 \quad \text{on } \{0,\pi\}\times[0,T].
\end{equation}
Throughout, $u_0\in H^1(D)$ is deterministic, $\|\cdot\|$ and $\langle\cdot,\cdot\rangle$ denote the norm and scalar product of $L^2(D)$, and $\|\varphi(t)\|_\infty=\sup_{x\in D}|\varphi(t,x)|$. Constants $C$ may change from line to line and depend only on $T$, $\varepsilon$, $\sigma$ unless stated otherwise.

\subsection{Spectral setting}

Set $e_0(x)=\pi^{-1/2}$ and, for $k\geq1$, $e_k(x)=(2/\pi)^{1/2}\cos(kx)$ and $g_k(x)=(2/\pi)^{1/2}\sin(kx)$. Both $(e_k)_{k\geq0}$ and $(g_k)_{k\geq1}$ are orthonormal bases of $L^2(D)$. Each $e_k$ satisfies \eqref{eq:bc},
\begin{equation}\label{eq:basis}
-\Delta e_k=k^2e_k,\qquad e_k'=-k\,g_k,\qquad g_k'=k\,e_k,\qquad |e_k|,|g_k|\leq(2/\pi)^{1/2},
\end{equation}
and $g_k(0)=g_k(\pi)=0$. For $w\in L^2(D)$ we write $\hat w_k=\langle w,e_k\rangle$. For $s\geq0$, the norm $\|w\|_{H^s}^2:=\sum_{k\geq0}(1+k^2)^s|\hat w_k|^2$ is equivalent to the usual Sobolev norm of $H^s(D)$ for $s\in\{0,1\}$; for larger $s$ it defines the domain of the corresponding power of $1-\Delta$ with the boundary conditions \eqref{eq:bc}. We only use the spectral norms below.

We define $A$ as the self-adjoint operator on $L^2(D)$ given by $Aw=-\sum_k\lambda_k\hat w_ke_k$, with domain $\mathcal D(A)=\{w\in L^2(D):\sum_k\lambda_k^2|\hat w_k|^2<\infty\}$, where $\lambda_k$ is given below. For $w\in C^6(D)$ satisfying \eqref{eq:bc}, three successive pairs of integrations by parts show that $\langle\Delta(1+\Delta)^2w,e_k\rangle=\langle w,\Delta(1+\Delta)^2e_k\rangle$, so that $A$ is the realization of the differential operator $\Delta(1+\Delta)^2$ with the boundary conditions \eqref{eq:bc}; this is the only realization used in the paper. It is diagonal in $(e_k)$:
\begin{equation}\label{eq:eigen}
A e_k = -\lambda_k e_k, \qquad \lambda_k := k^2(1-k^2)^2 \geq 0 .
\end{equation}
We also set $\|A^{1/2}w\|^2:=\sum_{k\geq0}\lambda_k|\hat w_k|^2$. For $w$ smooth and satisfying \eqref{eq:bc}, $\|A^{1/2}w\|^2=\|\partial_x(1+\Delta)w\|^2=-\langle Aw,w\rangle$, by \eqref{eq:basis}.

\begin{lemma}\label{lem:spectrum}
$\lambda_0=\lambda_1=0$, and $\lambda_k\geq \frac{9}{16}\,k^6$ for every $k\geq2$, with equality for $k=2$.
\end{lemma}

\begin{proof}
The first claim is immediate. For $k\geq2$, $k^2-1\geq\frac34k^2$, with equality for $k=2$, hence $\lambda_k=k^2(k^2-1)^2\geq\frac{9}{16}k^6$.
\end{proof}

\begin{remark}[Two neutral modes]\label{rem:neutral}
For the Cahn-Hilliard operator $-\Delta^2$ the only zero eigenvalue is the one of the constant mode. Here the mode $e_1(x)=(2/\pi)^{1/2}\cos x$ is also neutral. The constant mode is fixed by the conservation of the mean (Proposition~\ref{prop:mass}). The mode $e_1$ is neither damped by the linear part nor fixed by a conservation law; the dissipation $\|A^{1/2}w\|^2$ does not see it, and it requires the separate argument of Lemma~\ref{lem:interp}.
\end{remark}

\subsection{The Green function}

For $t>0$ and $x,y\in D$ let
\begin{equation}\label{eq:green}
\begin{split}
G(t,x,y)&=\sum_{k\geq0} e^{-\lambda_k t}\,e_k(x)e_k(y)\\
&= G^0(x,y)+\widetilde G(t,x,y),
\end{split}
\end{equation}
where $G^0(x,y)=\frac1\pi+\frac2\pi\cos x\cos y$ collects the neutral modes and
\[
\widetilde G(t,x,y)=\sum_{k\geq2}e^{-\lambda_kt}e_k(x)e_k(y).
\]
For $N\geq1$ we also use the truncated kernels $G_N(t,x,y)=\sum_{k=0}^Ne^{-\lambda_kt}e_k(x)e_k(y)$, and $S(t)w=\int_DG(t,\cdot,y)w(y)\,dy=\sum_ke^{-\lambda_kt}\hat w_ke_k$. By Lemma~\ref{lem:sums}(a) below, all the series of derivatives of $\widetilde G$ converge absolutely and uniformly on $[t_0,\infty)\times D\times D$ for every $t_0>0$, so that $G$ is smooth on $(0,\infty)\times D\times D$ and may be differentiated term by term.

\begin{lemma}[Pointwise bounds]\label{lem:green}
For all integers $i,j\geq0$ there exists $C_{i,j}=C_{i,j}(T)$ such that, for all $N\geq1$, $t\in(0,T]$ and $x,y\in D$,
\begin{equation}\label{eq:greenbound}
\big|\partial_x^i\partial_y^j \widetilde G(t,x,y)\big| \leq C_{i,j}\, t^{-(i+j+1)/6},\qquad
\big|\partial_x^i\partial_y^j G_N(t,x,y)\big| \leq C_{i,j}\big(1+t^{-(i+j+1)/6}\big).
\end{equation}
The same bound as for $G_N$ holds for $G$.
\end{lemma}

\begin{proof}
Fix $i,j\geq0$, set $m=i+j$ and let $t\in(0,T]$, $x,y\in D$. By \eqref{eq:basis}, every derivative of $e_k$ of order $n$ is, up to sign, $k^n$ times $e_k$ or $g_k$, so that $|\partial^ne_k|\leq(2/\pi)^{1/2}k^n$. Differentiating \eqref{eq:green} term by term (which is justified by the uniform convergence of the resulting series on $[t_0,T]\times D\times D$ for every $t_0>0$, see below) and using Lemma~\ref{lem:spectrum}, we obtain
\begin{align*}
\big|\partial_x^i\partial_y^j \widetilde G(t,x,y)\big|
&=\Big|\sum_{k\geq2}e^{-\lambda_kt}\,\partial_x^ie_k(x)\,\partial_y^je_k(y)\Big|
\leq\sum_{k\geq2}e^{-\lambda_kt}\,|\partial_x^ie_k(x)|\,|\partial_y^je_k(y)|\\
&\leq\frac2\pi\sum_{k\geq2}k^{m}\,e^{-\lambda_kt}
\leq\frac2\pi\sum_{k\geq2}k^{m}\,e^{-\frac9{16}k^6t}\\
&\leq\frac2\pi\Big(\int_0^\infty r^me^{-\frac9{16}r^6t}\,dr+\sup_{r\geq0}r^me^{-\frac9{16}r^6t}\Big)\\
&=\frac2\pi\Big(t^{-\frac{m+1}6}\int_0^\infty\rho^me^{-\frac9{16}\rho^6}\,d\rho+t^{-\frac m6}\sup_{\rho\geq0}\rho^me^{-\frac9{16}\rho^6}\Big)\\
&\leq\frac2\pi\Big(c_m+T^{1/6}c'_m\Big)\,t^{-\frac{m+1}6}=:C_{i,j}\,t^{-\frac{m+1}6},
\end{align*}
where the third inequality is Lemma~\ref{lem:sums}(a) (comparison of the sum of a unimodal function with its integral), the equality follows from the change of variables $r=t^{-1/6}\rho$, and the last inequality uses $t^{-m/6}=t^{1/6}t^{-(m+1)/6}\leq T^{1/6}t^{-(m+1)/6}$. The same computation, for $t\geq t_0$, bounds the series of absolute values by $\frac2\pi\sum_kk^me^{-\frac9{16}k^6t_0}<\infty$, which gives the uniform convergence used above.

For $G_N$ and $G$, the modes $k\geq2$ are bounded by the same sum of nonnegative terms, restricted to $k\leq N$ in the case of $G_N$. The neutral modes contribute $\partial_x^i\partial_y^j\big(\frac1\pi+e_1(x)e_1(y)\big)$, whose absolute value is at most $\frac1\pi+\frac2\pi\leq1$. Hence $|\partial_x^i\partial_y^jG_N(t,x,y)|\leq1+C_{i,j}t^{-(m+1)/6}\leq C_{i,j}(1+t^{-(m+1)/6})$ after enlarging $C_{i,j}$.
\end{proof}

\begin{remark}[Comparison with Cahn-Hilliard]
For $-\Delta^2$, $\lambda_k=k^4$ and the analogous bound is $Ct^{-(i+j+1)/4}$. The nonlinear term of the mild formulation involves $\partial_y^2G$, for which Lemma~\ref{lem:green} gives $C(1+t^{-1/2})$, against $Ct^{-3/4}$ for Cahn-Hilliard.
\end{remark}

\subsection{Mild solutions}

For $\varphi\in L^\infty([0,T]\times D)$ we set
\begin{equation}\label{eq:N}
N_\varphi(t,x)=\int_0^t\!\!\int_D \Delta_yG(t-s,x,y)\,\varphi(s,y)\,dy\,ds,
\end{equation}
where the integral converges absolutely by Lemma~\ref{lem:green} with $(i,j)=(0,2)$. The stochastic convolution is the Walsh integral
\begin{equation}\label{eq:Z}
Z(t,x)=-\sigma\int_0^t\!\!\int_D \partial_y G(t-s,x,y)\,W(ds,dy),
\end{equation}
which is well defined since $\int_0^t\int_D|\partial_yG(r,x,y)|^2dy\,dr<\infty$ (Lemma~\ref{lem:Zmoments}).

\begin{definition}\label{def:mild}
Let $\theta$ be a stopping time with values in $(0,T]$. A mild solution of \eqref{eq:pfc}-\eqref{eq:bc} on $[0,\theta)$ is an adapted random field $u$ such that, almost surely, $u\in C([0,t]\times D)$ for every $t<\theta$ and
\begin{equation}\label{eq:mild}
u(t,x)=S(t)u_0(x)+N_{f(u)}(t,x)+Z(t,x)\qquad\text{for all }t<\theta,\ x\in D.
\end{equation}
A mild solution on $[0,T]$ is defined in the same way with $\theta=T$ and $t\leq T$.
\end{definition}

\begin{remark}
Formula \eqref{eq:mild} is obtained from \eqref{eq:pfc} by Duhamel's principle and two integrations by parts in $y$: $\int_DG\,\Delta f(u)\,dy=\int_D\Delta_yG\,f(u)\,dy$ and $\int_DG\,\partial_y\dot W\,dy=-\int_D\partial_yG\,\dot W\,dy$. The boundary terms vanish because $\partial_yG(t,x,\cdot)$ vanishes at $0$ and $\pi$ (it is a series in the $g_k$) and because $\partial_xu=0$ there. We take \eqref{eq:mild} as the definition of a solution.
\end{remark}

\begin{proposition}[Conservation of the mean]\label{prop:mass}
If $u$ is a mild solution on $[0,\theta)$, then almost surely $\int_Du(t,x)\,dx=\int_Du_0(x)\,dx$ for all $t<\theta$.
\end{proposition}

\begin{proof}
Since $\int_De_k(x)\,dx=0$ for $k\geq1$, we have $\int_DS(t)u_0(x)\,dx=\pi^{1/2}\hat u_{0,0}=\int_Du_0$. By term by term integration, $\int_D\Delta_yG(r,x,y)\,dx=-\sum_{k\geq1}k^2e^{-\lambda_kr}e_k(y)\int_De_k(x)\,dx=0$ for $r>0$, so Fubini's theorem gives $\int_DN_\varphi(t,x)\,dx=0$. Likewise $\int_D\partial_yG(r,x,y)\,dx=0$, and the stochastic Fubini theorem \cite[Theorem~2.6]{Wa86} gives $\int_DZ(t,x)\,dx=0$ almost surely for each $t$; by continuity of $Z$ (Lemma~\ref{lem:Zholder}) this holds for all $t$ simultaneously.
\end{proof}

\subsection{The stochastic convolution}

For $k\geq1$ let $B_k(t)=\int_0^t\!\int_Dg_k(y)\,W(ds,dy)$. Since $(g_k)$ is orthonormal, the $B_k$ are independent standard $(\mathcal F_t)$-Brownian motions. Let
\[
z_1(t)=B_1(t),\qquad z_k(t)=\int_0^te^{-\lambda_k(t-s)}\,dB_k(s)\quad(k\geq2).
\]

\begin{lemma}[Spectral representation and moments]\label{lem:Zmoments}
For every $(t,x)$,
\begin{equation}\label{eq:Zseries}
Z(t,x)=\sigma\sum_{k\geq1}k\,e_k(x)\,z_k(t)\quad\text{in }L^2(\Omega),
\end{equation}
and
\begin{equation}\label{eq:normG}
\E|Z(t,x)|^2=\sigma^2\int_0^t\!\!\int_D|\partial_yG(r,x,y)|^2\,dy\,dr=\sigma^2\Big(e_1(x)^2\,t+\sum_{k\geq2}k^2e_k(x)^2\,\frac{1-e^{-2\lambda_kt}}{2\lambda_k}\Big)\leq C.
\end{equation}
\end{lemma}

\begin{proof}
By \eqref{eq:basis} and term by term differentiation, for $r>0$,
\[
\partial_yG(r,x,y)=-\sum_{k\geq1}k\,e^{-\lambda_kr}e_k(x)g_k(y),
\]
and by orthonormality of $(g_k)$,
\[
\int_D|\partial_yG(r,x,y)|^2\,dy=\sum_{k\geq1}k^2e^{-2\lambda_kr}e_k(x)^2 .
\]
Integrating in $r\in[0,t]$ gives the series in \eqref{eq:normG}, which is bounded by $\frac2\pi\big(T+\sum_{k\geq2}k^2/(2\lambda_k)\big)\leq C$ since $k^2/\lambda_k\leq\frac{16}9k^{-4}$. Hence $\partial_yG(t-\cdot,x,\cdot)\mathbf 1_{[0,t]}\in\cH$, and the partial sums of the series $-\sum_{k}k\,e^{-\lambda_k(t-s)}e_k(x)g_k(y)$ converge to it in $\cH$. Since the integral of the $n$-th partial sum is
\[
-\sigma\int_0^t\!\!\int_D\Big(-\sum_{k=1}^nk\,e^{-\lambda_k(t-s)}e_k(x)g_k(y)\Big)W(ds,dy)=\sigma\sum_{k=1}^nk\,e_k(x)\,z_k(t)
\]
(recall $\lambda_1=0$, so that the term $k=1$ is $\sigma e_1(x)B_1(t)$), the isometry of the Walsh integral shows that these partial sums converge in $L^2(\Omega)$ to $Z(t,x)$; this is \eqref{eq:Zseries}. The variance formula follows from the isometry.
\end{proof}

\begin{lemma}[H\"older continuity of $Z$]\label{lem:Zholder}
For all $s,t\in[0,T]$ and $x,y\in D$,
\[
\E|Z(t,x)-Z(t,y)|^2\leq C|x-y|^2,\qquad \E|Z(t,x)-Z(s,x)|^2\leq C|t-s|^{1/2}.
\]
Consequently, $Z$ has a version which is continuous on $[0,T]\times D$, $\gamma_1$-H\"older in $t$ and $\gamma_2$-H\"older in $x$ for all $\gamma_1<\frac14$, $\gamma_2<1$. We always use this version.
\end{lemma}

\begin{proof}
\emph{Space.} Let $t\in[0,T]$, $x,y\in D$ and $h=|x-y|$. By \eqref{eq:Zseries} and the independence of the centred variables $z_k(t)$,
\[
\E|Z(t,x)-Z(t,y)|^2=\sigma^2\sum_{k\geq1}k^2\,|e_k(x)-e_k(y)|^2\,\E|z_k(t)|^2 .
\]
By the mean value theorem, $|e_k(x)-e_k(y)|\leq(2/\pi)^{1/2}\min(kh,2)$. Also $\E|z_1(t)|^2=t$ and, for $k\geq2$, $\E|z_k(t)|^2=\frac{1-e^{-2\lambda_kt}}{2\lambda_k}\leq\frac1{2\lambda_k}\leq\frac89k^{-6}$ by Lemma~\ref{lem:spectrum}. Hence
\begin{align*}
\E|Z(t,x)-Z(t,y)|^2
&=\sigma^2\Big(|e_1(x)-e_1(y)|^2\,t+\sum_{k\geq2}k^2|e_k(x)-e_k(y)|^2\,\frac{1-e^{-2\lambda_kt}}{2\lambda_k}\Big)\\
&\leq\frac{2\sigma^2}\pi\Big(T\,h^2+\sum_{k\geq2}k^2\min(k^2h^2,4)\,\frac1{2\lambda_k}\Big)\\
&\leq\frac{2\sigma^2}\pi\Big(T\,h^2+\frac89\sum_{k\geq2}k^{-4}\min(k^2h^2,4)\Big)\\
&\leq\frac{2\sigma^2}\pi\Big(T+\frac89C_4\Big)h^2=:C\,h^2,
\end{align*}
where $C_4$ is the constant of Lemma~\ref{lem:sums}(c) with $p=4$, for which $\min(p-1,2)=2$.

\emph{Time.} Let $0\leq s<t\leq T$ and $h=t-s$. For $k\geq2$,
\[
z_k(t)-z_k(s)=\big(e^{-\lambda_kh}-1\big)z_k(s)+\int_s^te^{-\lambda_k(t-r)}\,dB_k(r),
\]
where the two terms are independent, the first being $\mathcal F_s$-measurable and the second independent of $\mathcal F_s$. Using the It\^o isometry and the elementary inequalities $(1-e^{-a})^2\leq1-e^{-a}\leq\min(a,1)$ for $a\geq0$,
\begin{align*}
\E|z_k(t)-z_k(s)|^2
&=\big(1-e^{-\lambda_kh}\big)^2\,\E|z_k(s)|^2+\int_s^te^{-2\lambda_k(t-r)}\,dr\\
&=\big(1-e^{-\lambda_kh}\big)^2\,\frac{1-e^{-2\lambda_ks}}{2\lambda_k}+\frac{1-e^{-2\lambda_kh}}{2\lambda_k}\\
&\leq\frac{1-e^{-\lambda_kh}}{2\lambda_k}+\frac{1-e^{-2\lambda_kh}}{2\lambda_k}
\leq\frac{\min(\lambda_kh,1)}{2\lambda_k}+\frac{\min(2\lambda_kh,1)}{2\lambda_k}\\
&\leq\frac12\min\Big(h,\frac1{\lambda_k}\Big)+\min\Big(h,\frac1{2\lambda_k}\Big)\leq2\min\Big(h,\frac1{\lambda_k}\Big).
\end{align*}
For $k=1$, $\E|z_1(t)-z_1(s)|^2=\E|B_1(t)-B_1(s)|^2=h$. Therefore, by independence and Lemma~\ref{lem:spectrum},
\begin{align*}
\E|Z(t,x)-Z(s,x)|^2
&=\sigma^2\sum_{k\geq1}k^2e_k(x)^2\,\E|z_k(t)-z_k(s)|^2\\
&\leq\frac{2\sigma^2}\pi\Big(h+2\sum_{k\geq2}k^2\min\Big(h,\frac1{\lambda_k}\Big)\Big)
\leq\frac{2\sigma^2}\pi\Big(h+2\sum_{k\geq2}k^2\min\Big(h,\frac{16}9k^{-6}\Big)\Big)\\
&\leq\frac{2\sigma^2}\pi\Big(h+\frac{32}9\sum_{k\geq2}k^2\min\big(h,k^{-6}\big)\Big)\\
&\leq\frac{2\sigma^2}\pi\Big(T^{1/2}h^{1/2}+\frac{32}9C'_2h^{1/2}\Big)=:C\,h^{1/2},
\end{align*}
where $C'_2$ is the constant of Lemma~\ref{lem:sums}(b) with $m=2$, for which $\frac{5-m}6=\frac12$, and we used $h=h^{1/2}h^{1/2}\leq T^{1/2}h^{1/2}$.

\emph{Continuity.} For $(t,x),(s,y)\in[0,T]\times D$, the triangle inequality and the two bounds above give
\[
\E|Z(t,x)-Z(s,y)|^2\leq2\E|Z(t,x)-Z(t,y)|^2+2\E|Z(t,y)-Z(s,y)|^2\leq C\big(|t-s|^{1/2}+|x-y|^2\big).
\]
Lemma~\ref{lem:kolmogorov} with $\alpha_1=\frac14$ and $\alpha_2=1$ gives the claim.
\end{proof}

\begin{lemma}[Spatial derivative of $Z$]\label{lem:Zx}
Almost surely, the series \eqref{eq:Zseries} and
\begin{equation}\label{eq:Zxseries}
Y(t,x):=-\sigma\sum_{k\geq1}k^2\,g_k(x)\,z_k(t)
\end{equation}
converge uniformly on $[0,T]\times D$, their sums are continuous, $Z(t,\cdot)\in C^1(D)$ and $\partial_xZ=Y$. In addition $\partial_xZ(t,0)=\partial_xZ(t,\pi)=0$, and
\begin{equation}\label{eq:KZ}
K_Z:=\sup_{t\in[0,T]}\big(\|Z(t)\|_\infty+\|\partial_xZ(t)\|_\infty\big)<\infty\quad\text{almost surely},\qquad \E\big[K_Z^p\big]<\infty\ \text{ for all }p\geq1.
\end{equation}
\end{lemma}

\begin{proof}
\emph{Step 1: moments and increments of $Y$.} By independence of the $z_k(t)$ and $|g_k|^2\leq\frac2\pi$,
\begin{align*}
\E|Y(t,x)|^2&=\sigma^2\sum_{k\geq1}k^4g_k(x)^2\,\E|z_k(t)|^2
\leq\frac{2\sigma^2}\pi\Big(T+\sum_{k\geq2}\frac{k^4}{2\lambda_k}\Big)
\leq\frac{2\sigma^2}\pi\Big(T+\frac89\sum_{k\geq2}k^{-2}\Big)<\infty,
\end{align*}
so the series \eqref{eq:Zxseries} converges in $L^2(\Omega)$ for each $(t,x)$. For the spatial increments, let $h=|x-y|$; since $|g_k(x)-g_k(y)|\leq(2/\pi)^{1/2}\min(kh,2)$,
\begin{align*}
\E|Y(t,x)-Y(t,y)|^2
&=\sigma^2\sum_{k\geq1}k^4|g_k(x)-g_k(y)|^2\,\E|z_k(t)|^2\\
&\leq\frac{2\sigma^2}\pi\Big(T\,h^2+\sum_{k\geq2}k^4\min(k^2h^2,4)\frac1{2\lambda_k}\Big)\\
&\leq\frac{2\sigma^2}\pi\Big(T\,h^2+\frac89\sum_{k\geq2}k^{-2}\min(k^2h^2,4)\Big)\\
&\leq\frac{2\sigma^2}\pi\Big(T\pi+\frac89C_2\Big)h ,
\end{align*}
by Lemma~\ref{lem:sums}(c) with $p=2$ and $h^2\leq\pi h$. For the time increments, with $h=|t-s|$ and the bound $\E|z_k(t)-z_k(s)|^2\leq2\min(h,\lambda_k^{-1})$ from the proof of Lemma~\ref{lem:Zholder},
\begin{align*}
\E|Y(t,x)-Y(s,x)|^2
&\leq\frac{2\sigma^2}\pi\Big(h+2\sum_{k\geq2}k^4\min\Big(h,\frac{16}9k^{-6}\Big)\Big)
\leq\frac{2\sigma^2}\pi\Big(h+\frac{32}9\sum_{k\geq2}k^4\min(h,k^{-6})\Big)\\
&\leq\frac{2\sigma^2}\pi\Big(T^{5/6}h^{1/6}+\frac{32}9C'_4h^{1/6}\Big),
\end{align*}
by Lemma~\ref{lem:sums}(b) with $m=4$, for which $\frac{5-m}6=\frac16$. Since $Y$ is a centred Gaussian field, Lemma~\ref{lem:kolmogorov} with $\alpha_1=\frac1{12}$, $\alpha_2=\frac12$ gives a continuous version of $Y$.

\emph{Step 2: uniform convergence.} The partial sums $Z^n=\sigma\sum_{k\leq n}k\,e_kz_k$ and $Y^n=-\sigma\sum_{k\leq n}k^2g_kz_k$ are sums of independent, symmetric random variables with values in the separable Banach space $C([0,T]\times D)$, since each $z_k$ has continuous paths. The dual of $C([0,T]\times D)$ is the space of finite signed measures $\nu$ on $[0,T]\times D$. For such $\nu$, by the Cauchy-Schwarz inequality and Fubini's theorem, using the continuous versions of $Z$ and $Y$,
\begin{align*}
\E\Big|\int(Z^n-Z)\,d\nu\Big|^2
&\leq|\nu|([0,T]\times D)\int\E|Z^n(t,x)-Z(t,x)|^2\,d|\nu|(t,x)\\
&=|\nu|([0,T]\times D)\int\sigma^2\sum_{k>n}k^2e_k(x)^2\,\E|z_k(t)|^2\,d|\nu|(t,x)\\
&\leq\frac{16\sigma^2}{9\pi}\,|\nu|([0,T]\times D)^2\sum_{k>n}k^{-4},
\end{align*}
which tends to $0$; the same computation for $Y$ gives the bound $\frac{16\sigma^2}{9\pi}|\nu|([0,T]\times D)^2\sum_{k>n}k^{-2}$. Hence $\int Z^nd\nu\to\int Z\,d\nu$ and $\int Y^nd\nu\to\int Y\,d\nu$ in probability for every $\nu$, and the It\^o-Nisio theorem (Lemma~\ref{lem:itonisio}) shows that, almost surely, $Z^n\to Z$ and $Y^n\to Y$ uniformly on $[0,T]\times D$.

\emph{Step 3: derivative.} By \eqref{eq:basis}, $\partial_xZ^n=\sigma\sum_{k\leq n}k\,e_k'z_k=-\sigma\sum_{k\leq n}k^2g_kz_k=Y^n$. For fixed $t$, $Z^n(t,\cdot)\to Z(t,\cdot)$ and $\partial_xZ^n(t,\cdot)\to Y(t,\cdot)$ uniformly on $D$; by the classical theorem on uniform convergence of derivatives, $Z(t,\cdot)\in C^1(D)$ and $\partial_xZ(t,\cdot)=Y(t,\cdot)$. Since $g_k(0)=g_k(\pi)=0$, $Y^n(t,0)=Y^n(t,\pi)=0$ for all $n$, and the uniform limit also vanishes at the boundary.

\emph{Step 4: Gaussian integrability of $K_Z$.} Let $B=C([0,T]\times D)^2$ with the norm $\|(a,b)\|_B=\|a\|_\infty+\|b\|_\infty$; it is a separable Banach space. By Steps~2 and~3, $X:=(Z,Y)$ is a centred Gaussian random variable with values in $B$ (it is the almost sure limit in $B$ of the Gaussian vectors $(Z^n,Y^n)$), and
\[
K_Z=\sup_{t\in[0,T]}\big(\|Z(t)\|_\infty+\|\partial_xZ(t)\|_\infty\big)\leq\|Z\|_\infty+\|Y\|_\infty=\|X\|_B .
\]
Let $s(X)=\sup\{(\E\langle\phi,X\rangle^2)^{1/2}:\phi\in B^*,\ \|\phi\|_{B^*}\leq1\}$, which is finite since $\E\|X\|_B^2<\infty$. By Fernique's theorem in the form of \cite[Corollary~3.2]{LT91},
\[
\E\exp\Big(\frac{\|X\|_B^2}{2\alpha^2}\Big)<\infty\qquad\text{for every }\alpha>s(X).
\]
Hence there is $c>0$, for instance $c=1/(8s(X)^2)$, such that $\E\exp(cK_Z^2)\leq\E\exp(c\|X\|_B^2)<\infty$. In particular, since $r^p\leq C_{p,c}e^{cr^2}$ for $r\geq0$, $\E[K_Z^p]<\infty$ for every $p\geq1$.
\end{proof}

\begin{remark}
Lemma~\ref{lem:Zx} improves the spatial part of Lemma~\ref{lem:Zholder}: because of the conservative form of the noise, $Z(t,\cdot)$ is $C^1$, and $\partial_xZ$ is $\gamma$-H\"older in $x$ for every $\gamma<\frac12$ and in $t$ for every $\gamma<\frac1{12}$. This is the key property used in Section~\ref{sec:global}.
\end{remark}
\subsection{Mild and weak solutions}

We now show that mild solutions coincide with weak solutions in the usual sense. Let $\mathcal D$ be the set of $\phi\in C^6(D)$ such that $\phi'=\phi'''=\phi^{(5)}=0$ at $0$ and $\pi$. For $\phi\in\mathcal D$, two integrations by parts in each of the three factors of $A=\Delta(1+\Delta)^2$ give $\langle\phi,e_k\rangle$-coefficients
\begin{equation}\label{eq:coefD}
\langle A\phi,e_k\rangle=\langle\phi,Ae_k\rangle=-\lambda_k\hat\phi_k,\qquad \langle\Delta\phi,e_k\rangle=-k^2\hat\phi_k,\qquad \langle\phi',g_k\rangle=-k\hat\phi_k ,
\end{equation}
where the boundary terms vanish thanks to the conditions on $\phi$ and to $e_k'(0)=e_k'(\pi)=0$, $g_k(0)=g_k(\pi)=0$.

\begin{proposition}[Weak formulation]\label{prop:weak}
Let $\theta$ be a stopping time with values in $(0,T]$ and $u$ an adapted random field with continuous paths on $[0,t]\times D$ for every $t<\theta$. Then $u$ is a mild solution on $[0,\theta)$ if and only if, almost surely, for every $\phi\in\mathcal D$ and every $t<\theta$,
\begin{equation}\label{eq:weak}
\langle u(t),\phi\rangle=\langle u_0,\phi\rangle+\int_0^t\langle u(s),A\phi\rangle\,ds+\int_0^t\langle f(u(s)),\Delta\phi\rangle\,ds-\sigma\int_0^t\!\!\int_D\phi'(y)\,W(ds,dy).
\end{equation}
\end{proposition}

\begin{proof}
\emph{Step 1: coefficients of the three terms of \eqref{eq:mild}.} Let $\varphi\in L^\infty([0,T]\times D)$ and $k\geq0$. By the orthonormality of $(e_k)$, $\langle S(t)u_0,e_k\rangle=e^{-\lambda_kt}\hat u_{0,k}$. Since $\int_D\Delta_yG(r,x,y)e_k(x)\,dx=-k^2e^{-\lambda_kr}e_k(y)$ for $r>0$ (term by term integration of a uniformly convergent series), Fubini's theorem gives
\[
\langle N_\varphi(t),e_k\rangle=\int_0^t\!\!\int_D\Big(\int_D\Delta_yG(t-s,x,y)e_k(x)\,dx\Big)\varphi(s,y)\,dy\,ds=-k^2\int_0^te^{-\lambda_k(t-s)}\langle\varphi(s),e_k\rangle\,ds .
\]
Finally, by Lemma~\ref{lem:Zx}, the series \eqref{eq:Zseries} converges uniformly on $D$, so it can be integrated term by term against $e_k$: $\langle Z(t),e_k\rangle=\sigma k\,z_k(t)$ for $k\geq1$, and $\langle Z(t),e_0\rangle=0$.

\emph{Step 2: the coefficient equations.} Let $\xi_k(t)=\langle u(t),e_k\rangle$, $\varphi=f(u)$ and $\psi_k(s)=\langle\varphi(s),e_k\rangle$. By Step~1, $u$ is a mild solution if and only if, for every $k\geq0$ and $t<\theta$,
\begin{equation}\label{eq:varconst}
\xi_k(t)=e^{-\lambda_kt}\hat u_{0,k}-k^2\int_0^te^{-\lambda_k(t-s)}\psi_k(s)\,ds+\sigma k\,z_k(t),
\end{equation}
since a continuous function on $D$ is determined by its coefficients. On the other hand, by \eqref{eq:coefD} and the identity $\int_0^t\!\int_D\phi'(y)W(ds,dy)=\sum_{k\geq1}\langle\phi',g_k\rangle B_k(t)=-\sum_{k\geq1}k\hat\phi_kB_k(t)$, applied to $\phi=e_k$, \eqref{eq:weak} with $\phi=e_k$ reads
\begin{equation}\label{eq:weakk}
\xi_k(t)=\hat u_{0,k}-\lambda_k\int_0^t\xi_k(s)\,ds-k^2\int_0^t\psi_k(s)\,ds+\sigma k\,B_k(t).
\end{equation}

\emph{Step 3: \eqref{eq:varconst} is equivalent to \eqref{eq:weakk}.} Assume \eqref{eq:varconst}. Recall that $z_k(t)=B_k(t)-\lambda_k\int_0^tz_k(s)\,ds$ (the integral form of the Ornstein-Uhlenbeck equation, valid also for $k=1$ where $\lambda_1=0$). Integrating \eqref{eq:varconst} in time and using Fubini's theorem,
\begin{align*}
\lambda_k\int_0^t\xi_k(s)\,ds
&=\lambda_k\int_0^te^{-\lambda_ks}\,ds\;\hat u_{0,k}-k^2\int_0^t\Big(\lambda_k\int_r^te^{-\lambda_k(s-r)}ds\Big)\psi_k(r)\,dr+\sigma k\lambda_k\int_0^tz_k(s)\,ds\\
&=\big(1-e^{-\lambda_kt}\big)\hat u_{0,k}-k^2\int_0^t\big(1-e^{-\lambda_k(t-r)}\big)\psi_k(r)\,dr+\sigma k\big(B_k(t)-z_k(t)\big)\\
&=\hat u_{0,k}-k^2\int_0^t\psi_k(r)\,dr+\sigma kB_k(t)\\
&\qquad-\Big(e^{-\lambda_kt}\hat u_{0,k}-k^2\int_0^te^{-\lambda_k(t-r)}\psi_k(r)\,dr+\sigma kz_k(t)\Big)\\
&=\hat u_{0,k}-k^2\int_0^t\psi_k(r)\,dr+\sigma kB_k(t)-\xi_k(t),
\end{align*}
which is \eqref{eq:weakk}. Conversely, if \eqref{eq:weakk} holds, then $\eta_k:=\xi_k-\sigma kz_k$ is continuous and satisfies $\eta_k(t)=\hat u_{0,k}-\lambda_k\int_0^t\eta_k(s)\,ds-k^2\int_0^t\psi_k(s)\,ds$, that is $\eta_k'=-\lambda_k\eta_k-k^2\psi_k$ in the sense of absolutely continuous functions, with $\psi_k$ continuous; the variation of constants formula gives $\eta_k(t)=e^{-\lambda_kt}\hat u_{0,k}-k^2\int_0^te^{-\lambda_k(t-s)}\psi_k(s)\,ds$, which is \eqref{eq:varconst}.

\emph{Step 4: from $e_k$ to $\phi\in\mathcal D$, with a single exceptional set.} We first give a pathwise meaning to the stochastic term of \eqref{eq:weak}, simultaneously for all $\phi\in\mathcal D$. By the Borel-Cantelli lemma and the Gaussian tail of $\sup_{t\leq T}|B_k(t)|$, there is an event $\Omega_1$ of probability one and a random constant $C_1(\omega)$ such that $\sup_{t\leq T}|B_k(t)|\leq C_1(\omega)\big(\log(k+1)\big)^{1/2}$ for all $k\geq1$ on $\Omega_1$. For $\phi\in\mathcal D$, the sequence $a_k:=\lambda_k\hat\phi_k=-\langle A\phi,e_k\rangle$ is square summable, so that, by Lemma~\ref{lem:spectrum} and the Cauchy-Schwarz inequality, on $\Omega_1$,
\[
\sum_{k\geq2}k\,|\hat\phi_k|\sup_{t\leq T}|B_k(t)|\leq\frac{16}9\,C_1(\omega)\sum_{k\geq2}|a_k|\,k^{-5}\big(\log(k+1)\big)^{1/2}\leq C\,C_1(\omega)\,\|A\phi\|<\infty .
\]
Hence the series $I_\phi(t):=\sum_{k\geq1}k\hat\phi_kB_k(t)$ converges absolutely and uniformly on $[0,T]$, for every $\phi\in\mathcal D$, on the single event $\Omega_1$; $I_\phi$ is continuous in $t$, linear in $\phi$, and, by \eqref{eq:coefD} and the isometry, $I_\phi(t)=-\int_0^t\!\int_D\phi'(y)\,W(ds,dy)$ almost surely for each $t$. We use $\sigma I_\phi(t)$ as the stochastic term of \eqref{eq:weak}.

Let $\Omega_2\subset\Omega_1$ be the event of probability one on which \eqref{eq:weakk} holds for all $k\geq0$ and all $t<\theta$ (a countable family of identities between continuous functions of $t$). On $\Omega_2$, fix $\phi\in\mathcal D$ and $t<\theta$, multiply \eqref{eq:weakk} by $\hat\phi_k$ and sum over $k$. By Parseval's identity and \eqref{eq:coefD},
\[
\sum_k\hat\phi_k\xi_k(t)=\langle u(t),\phi\rangle,\qquad
-\sum_k\lambda_k\hat\phi_k\xi_k(s)=\langle u(s),A\phi\rangle,\qquad
-\sum_kk^2\hat\phi_k\psi_k(s)=\langle f(u(s)),\Delta\phi\rangle,
\]
and the last two series are bounded, uniformly in $s\in[0,t]$, by $\|u(s)\|\,\|A\phi\|$ and $\|f(u(s))\|\,\|\Delta\phi\|$ through the Cauchy-Schwarz inequality; by dominated convergence, the sums and the time integrals can be interchanged. The noise terms sum to $\sigma I_\phi(t)$. This gives \eqref{eq:weak} for every $\phi\in\mathcal D$ and $t<\theta$ on $\Omega_2$. The converse direction only uses $\phi=e_k$, which belongs to $\mathcal D$.
\end{proof}

\section{Local existence and uniqueness}\label{sec:local}

Since the noise is additive, \eqref{eq:mild} can be solved pathwise. We fix $\omega$ in the event of probability one on which the conclusions of Lemmas~\ref{lem:Zholder} and~\ref{lem:Zx} hold, so that $Z\in C([0,T]\times D)$.

\subsection{The deterministic terms}

For $N\geq1$ we write $P_Nw=\sum_{k\leq N}\hat w_ke_k$ and, for $\varphi\in L^\infty([0,T]\times D)$,
\[
N^N_\varphi(t,x)=\int_0^t\!\!\int_D\Delta_yG_N(t-s,x,y)\,\varphi(s,y)\,dy\,ds .
\]

\begin{lemma}[Initial term]\label{lem:initial}
For $u_0\in H^1(D)$, the series $S(t)u_0(x)=\sum_{k\geq0}e^{-\lambda_kt}\hat u_{0,k}e_k(x)$ converges uniformly on $[0,T]\times D$. Consequently $S(\cdot)u_0\in C([0,T]\times D)$, and for all $N\geq1$,
\begin{equation}\label{eq:initial}
\begin{aligned}
&\sup_{t,x}|S(t)u_0(x)|+\sup_{t,x}|S(t)P_Nu_0(x)|\leq C\|u_0\|_{H^1},\\
&\sup_{t,x}|S(t)(u_0-P_Nu_0)(x)|\leq C\Big(\sum_{k>N}(1+k^2)|\hat u_{0,k}|^2\Big)^{1/2}.
\end{aligned}
\end{equation}
\end{lemma}

\begin{proof}
For $n\geq0$, since $e^{-\lambda_kt}\leq1$ and $|e_k|\leq(2/\pi)^{1/2}$, the Cauchy-Schwarz inequality gives
\[
\sum_{k>n}\big|e^{-\lambda_kt}\hat u_{0,k}e_k(x)\big|\leq C\Big(\sum_{k>n}(1+k^2)|\hat u_{0,k}|^2\Big)^{1/2}\Big(\sum_{k\geq1}(1+k^2)^{-1}\Big)^{1/2},
\]
uniformly in $(t,x)$. This tends to $0$ as $n\to\infty$, which gives the uniform convergence, the continuity of the sum, and the second bound in \eqref{eq:initial}. The first bound is the case $n=-1$, applied to $u_0$ and to $P_Nu_0$.
\end{proof}

\begin{lemma}[Deterministic convolution]\label{lem:conv}
There exists $C=C(T)$ such that, for all $\varphi\in L^\infty([0,T]\times D)$, $N\geq1$, $0\leq t\leq t+h\leq T$ and $x,x'\in D$, and for $M=N_\varphi$ as well as $M=N^N_\varphi$,
\begin{align}
|M(t,x)|&\leq C\int_0^t\big(1+(t-s)^{-1/2}\big)\,\|\varphi(s)\|_\infty\,ds, \label{eq:conv1}\\
|M(t,x)-M(t,x')|&\leq C\,\|\varphi\|_\infty\,|x-x'|, \label{eq:conv2}\\
|M(t+h,x)-M(t,x)|&\leq C\,\|\varphi\|_\infty\,h^{1/2}, \label{eq:conv3}
\end{align}
where $\|\varphi\|_\infty=\sup_{s,y}|\varphi(s,y)|$. In particular $N_\varphi,N^N_\varphi\in C([0,T]\times D)$. Moreover
\begin{equation}\label{eq:conv4}
\sup_{t,x}\big|N_\varphi(t,x)-N^N_\varphi(t,x)\big|\leq 2\,\|\varphi\|_\infty\sum_{k>N}\frac{k^2}{\lambda_k}\xrightarrow[N\to\infty]{}0 .
\end{equation}
\end{lemma}

\begin{proof}
We write the proof for $M=N_\varphi$. For $M=N^N_\varphi$ it is identical: all the sums over $k$ below are restricted to $k\leq N$, and all the bounds are obtained from sums of nonnegative terms, which can only decrease.

\emph{Proof of \eqref{eq:conv1}.} By Lemma~\ref{lem:green} with $(i,j)=(0,2)$, $|\Delta_yG(r,x,y)|\leq C_{0,2}(1+r^{-1/2})$. Hence
\begin{align*}
|N_\varphi(t,x)|&\leq\int_0^t\!\!\int_D|\Delta_yG(t-s,x,y)|\,|\varphi(s,y)|\,dy\,ds
\leq\int_0^t\Big(\int_D|\Delta_yG(t-s,x,y)|\,dy\Big)\|\varphi(s)\|_\infty\,ds\\
&\leq\pi C_{0,2}\int_0^t\big(1+(t-s)^{-1/2}\big)\|\varphi(s)\|_\infty\,ds .
\end{align*}

\emph{Proof of \eqref{eq:conv2}.} By the mean value theorem and Lemma~\ref{lem:green} with $(i,j)=(1,2)$,
\begin{align*}
|N_\varphi(t,x)-N_\varphi(t,x')|
&\leq\int_0^t\!\!\int_D\big|\Delta_yG(t-s,x,y)-\Delta_yG(t-s,x',y)\big|\,|\varphi(s,y)|\,dy\,ds\\
&\leq|x-x'|\,\|\varphi\|_\infty\int_0^t\!\!\int_D\sup_{\xi\in D}\big|\partial_x\Delta_yG(t-s,\xi,y)\big|\,dy\,ds\\
&\leq\pi C_{1,2}\,|x-x'|\,\|\varphi\|_\infty\int_0^T\big(1+r^{-2/3}\big)\,dr\\
&=\pi C_{1,2}\big(T+3T^{1/3}\big)\,\|\varphi\|_\infty\,|x-x'| .
\end{align*}

\emph{Proof of \eqref{eq:conv3}.} Write $M(t+h,x)-M(t,x)=I_1+I_2$, with
\begin{align*}
I_1&=\int_t^{t+h}\!\!\int_D\Delta_yG(t+h-s,x,y)\varphi(s,y)\,dy\,ds,\\
I_2&=\int_0^t\!\!\int_D\big[\Delta_yG(t+h-s,x,y)-\Delta_yG(t-s,x,y)\big]\varphi(s,y)\,dy\,ds .
\end{align*}
As in the proof of \eqref{eq:conv1},
\begin{align*}
|I_1|&\leq\pi C_{0,2}\|\varphi\|_\infty\int_t^{t+h}\big(1+(t+h-s)^{-1/2}\big)ds\\
&=\pi C_{0,2}\|\varphi\|_\infty\big(h+2h^{1/2}\big)\leq\pi C_{0,2}\big(T^{1/2}+2\big)\|\varphi\|_\infty h^{1/2}.
\end{align*}
For $I_2$, by \eqref{eq:basis}, $\Delta_yG(r,x,y)=-\sum_{k\geq1}k^2e^{-\lambda_kr}e_k(x)e_k(y)$ for $r>0$. The term $k=1$ does not depend on $r$ because $\lambda_1=0$, and it cancels in the difference. Since $e^{-\lambda_k(r+h)}\leq e^{-\lambda_kr}$, with the change of variables $r=t-s$,
\begin{align*}
|I_2|&=\Big|\int_0^t\!\!\int_D\sum_{k\geq2}k^2\big(e^{-\lambda_k(t-s)}-e^{-\lambda_k(t+h-s)}\big)e_k(x)e_k(y)\varphi(s,y)\,dy\,ds\Big|\\
&\leq\frac2\pi\,\pi\,\|\varphi\|_\infty\sum_{k\geq2}k^2\int_0^t\big(e^{-\lambda_kr}-e^{-\lambda_k(r+h)}\big)dr
\leq2\|\varphi\|_\infty\sum_{k\geq2}k^2\big(1-e^{-\lambda_kh}\big)\int_0^\infty e^{-\lambda_kr}\,dr\\
&=2\|\varphi\|_\infty\sum_{k\geq2}k^2\,\frac{1-e^{-\lambda_kh}}{\lambda_k}
\leq2\|\varphi\|_\infty\sum_{k\geq2}k^2\min\Big(h,\frac1{\lambda_k}\Big)\\
&\leq2\|\varphi\|_\infty\sum_{k\geq2}k^2\min\Big(h,\frac{16}9k^{-6}\Big)
\leq\frac{32}9\|\varphi\|_\infty\sum_{k\geq2}k^2\min\big(h,k^{-6}\big)
\leq\frac{32}9C'_2\,\|\varphi\|_\infty\,h^{1/2},
\end{align*}
where we used $1-e^{-a}\leq\min(a,1)$ and Lemma~\ref{lem:sums}(b) with $m=2$. The interchange of sum and integrals is justified by the absolute convergence of the series for $r>0$ and the integrability of its sum of absolute values, which is bounded by $C(1+r^{-1/2})$.

\emph{Continuity.} By \eqref{eq:conv2} and \eqref{eq:conv3}, $|M(t,x)-M(s,y)|\leq C\|\varphi\|_\infty(|t-s|^{1/2}+|x-y|)$, so $M\in C([0,T]\times D)$.

\emph{Proof of \eqref{eq:conv4}.} Since $\Delta_yG-\Delta_yG_N=-\sum_{k>N}k^2e^{-\lambda_kr}e_k(x)e_k(y)$,
\begin{align*}
|N_\varphi(t,x)-N^N_\varphi(t,x)|
&\leq\int_0^t\!\!\int_D\sum_{k>N}k^2e^{-\lambda_k(t-s)}|e_k(x)||e_k(y)||\varphi(s,y)|\,dy\,ds\\
&\leq\frac2\pi\,\pi\,\|\varphi\|_\infty\sum_{k>N}k^2\int_0^te^{-\lambda_kr}\,dr
\leq2\|\varphi\|_\infty\sum_{k>N}\frac{k^2}{\lambda_k}\leq\frac{32}9\|\varphi\|_\infty\sum_{k>N}k^{-4},
\end{align*}
which tends to $0$ as $N\to\infty$, uniformly in $(t,x)$.
\end{proof}

\subsection{Truncation}

Fix a nondecreasing function $\rho\in C^\infty(\R)$ with $\rho(u)=u$ for $|u|\leq1$, $\rho$ constant on $(-\infty,-2]$ and on $[2,\infty)$, and $0\leq\rho'\leq1$. For $R\geq1$ let $\rho_R(u)=R\rho(u/R)$ and
\[
f_R=f\circ\rho_R .
\]
Then $f_R\in C^\infty(\R)$, $f_R=f$ on $[-R,R]$, $|\rho_R|\leq2R$, and $f_R$, $f_R'$ and $f_R''$ are bounded. We set
\[
L_R:=\sup_{u\in\R}|f_R'(u)|\leq3(2R)^2+|\varepsilon|,\qquad B_R:=\sup_{u\in\R}|f_R(u)|\leq(2R)^3+2|\varepsilon|R .
\]

\subsection{The truncated equation}

\begin{proposition}\label{prop:trunc}
For every $R\geq1$, the equation
\begin{equation}\label{eq:trunc}
u_R(t,x)=S(t)u_0(x)+N_{f_R(u_R)}(t,x)+Z(t,x),\qquad (t,x)\in[0,T]\times D,
\end{equation}
has a unique solution $u_R\in C([0,T]\times D)$. The process $u_R$ is adapted, and the Picard iterates $u^0=0$, $u^{n+1}=S(\cdot)u_0+N_{f_R(u^n)}+Z$ satisfy
\begin{equation}\label{eq:picard}
\sup_{t,x}|u^n(t,x)-u_R(t,x)|\leq C_R\,2^{-n}\Big(\|u_0\|_{H^1}+B_R+\sup_{t,x}|Z(t,x)|\Big)
\end{equation}
with a deterministic constant $C_R$.
\end{proposition}

\begin{proof}
\emph{Step 1: $\Phi$ maps $E$ into $E$.} Let $E=C([0,T]\times D)$ and, for $v\in E$, $\Phi(v)=S(\cdot)u_0+N_{f_R(v)}+Z$. Since $f_R$ is continuous and bounded by $B_R$, $f_R(v)\in L^\infty([0,T]\times D)$. By Lemma~\ref{lem:initial}, $S(\cdot)u_0\in E$; by Lemma~\ref{lem:conv}, $N_{f_R(v)}\in E$; and $Z\in E$ for the fixed $\omega$. Hence $\Phi(v)\in E$.

\emph{Step 2: contraction.} For $\beta>0$ let $\|v\|_\beta=\sup_{t\in[0,T]}e^{-\beta t}\|v(t)\|_\infty$; since $e^{-\beta T}\|v\|_\infty\leq\|v\|_\beta\leq\|v\|_\infty$, this norm is equivalent to the sup norm and $(E,\|\cdot\|_\beta)$ is a Banach space. For $v,w\in E$, by linearity of $\varphi\mapsto N_\varphi$, \eqref{eq:conv1} and $|f_R(v)-f_R(w)|\leq L_R|v-w|$,
\begin{align*}
e^{-\beta t}\|\Phi(v)(t)-\Phi(w)(t)\|_\infty
&=e^{-\beta t}\big\|N_{f_R(v)-f_R(w)}(t)\big\|_\infty\\
&\leq Ce^{-\beta t}\int_0^t\big(1+(t-s)^{-1/2}\big)\,\|f_R(v(s))-f_R(w(s))\|_\infty\,ds\\
&\leq CL_R\int_0^t\big(1+(t-s)^{-1/2}\big)e^{-\beta(t-s)}\,e^{-\beta s}\|v(s)-w(s)\|_\infty\,ds\\
&\leq CL_R\,\|v-w\|_\beta\int_0^t\big(1+r^{-1/2}\big)e^{-\beta r}\,dr\\
&\leq CL_R\,\|v-w\|_\beta\Big(\int_0^\infty e^{-\beta r}\,dr+\int_0^\infty r^{-1/2}e^{-\beta r}\,dr\Big)\\
&=CL_R\big(\beta^{-1}+\Gamma(\tfrac12)\beta^{-1/2}\big)\|v-w\|_\beta .
\end{align*}
We choose $\beta=\beta_R$, deterministic, such that $CL_R\big(\beta_R^{-1}+\sqrt\pi\,\beta_R^{-1/2}\big)\leq\frac12$. Taking the supremum over $t$ gives $\|\Phi(v)-\Phi(w)\|_{\beta_R}\leq\frac12\|v-w\|_{\beta_R}$.

\emph{Step 3: fixed point and rate.} The Banach fixed point theorem gives a unique $u_R\in E$ with $\Phi(u_R)=u_R$, and the Picard iterates satisfy $\|u^{n+1}-u^n\|_{\beta_R}\leq2^{-n}\|u^1-u^0\|_{\beta_R}$, hence
\[
\|u^n-u_R\|_{\beta_R}\leq\sum_{m\geq n}\|u^{m+1}-u^m\|_{\beta_R}\leq2^{1-n}\|u^1-u^0\|_{\beta_R}=2^{1-n}\|u^1\|_{\beta_R}.
\]
Since $u^1=S(\cdot)u_0+N_{f_R(0)}+Z$, \eqref{eq:initial} and \eqref{eq:conv1} give
\[
\|u^1\|_{\beta_R}\leq\|u^1\|_\infty\leq C\|u_0\|_{H^1}+C\big(T+2T^{1/2}\big)B_R+\sup_{t,x}|Z(t,x)|,
\]
and $\|u^n-u_R\|_\infty\leq e^{\beta_RT}\|u^n-u_R\|_{\beta_R}$. This proves \eqref{eq:picard}.

\emph{Step 4: adaptedness.} $\Phi(v)(t,\cdot)$ only depends on $(v(s),Z(s))_{s\leq t}$. Since $u^0=0$ is adapted, by induction every $u^n(t,x)$ is $\mathcal F_t$-measurable, and so is the pointwise limit $u_R(t,x)$.
\end{proof}

\subsection{Local solution and blow-up alternative}

For $R\geq1$ let $\tau_R=\inf\{t\in[0,T]:\|u_R(t)\|_\infty\geq R\}$, with $\inf\emptyset=T$. Since $u_R$ is continuous and adapted, $\tau_R$ is a stopping time. Since $u_R(0)=u_0\in H^1(D)\subset C(D)$, we only consider $R>\|u_0\|_\infty$; then $\tau_R>0$ and, by continuity, $\|u_R(t)\|_\infty\leq R$ for all $t\in[0,\tau_R]$.

\begin{theorem}\label{thm:local}
Almost surely, the following holds.
\begin{enumerate}
\item For $\|u_0\|_\infty<R<R'$, $u_R=u_{R'}$ on $[0,\tau_R]\times D$, and $\tau_R\leq\tau_{R'}$, with $\tau_R<\tau_{R'}$ if $\tau_R<T$.
\item Let $\tau=\lim_{R\to\infty}\tau_R$, and $u=u_R$ on $[0,\tau_R]\times D$. Then $u$ is a mild solution on $[0,\tau)$ in the sense of Definition~\ref{def:mild}.
\item (Uniqueness) If $\tilde u$ is a mild solution on $[0,\theta)$, then $\tilde u=u$ on $[0,\theta\wedge\tau)\times D$.
\item (Blow-up alternative) On the event $\{\tau<T\}$, $\sup_{t<\tau}\|u(t)\|_\infty=+\infty$. Equivalently, if $\sup_{t<\tau}\|u(t)\|_\infty<\infty$, then $\tau=T$, $\tau_R=T$ for all large $R$, and $u$ is a mild solution on $[0,T]$.
\end{enumerate}
\end{theorem}

\begin{proof}
(1) For $t\leq\tau_R$ we have $\|u_R(t)\|_\infty\leq R$, hence $f_R(u_R)=f(u_R)=f_{R'}(u_R)$ on $[0,\tau_R]\times D$. Since $N_\varphi(t)$ only involves $\varphi$ on $[0,t]$, the difference $w=u_R-u_{R'}$ satisfies, for $t\leq\tau_R$,
\[
\|w(t)\|_\infty=\big\|N_{f_{R'}(u_R)-f_{R'}(u_{R'})}(t)\big\|_\infty\leq CL_{R'}\int_0^t\big(1+(t-s)^{-1/2}\big)\|w(s)\|_\infty\,ds .
\]
Lemma~\ref{lem:gronwall} with $a=0$ gives $w=0$ on $[0,\tau_R]$. Hence $\|u_{R'}(t)\|_\infty=\|u_R(t)\|_\infty<R<R'$ for $t<\tau_R$, and $\|u_{R'}(\tau_R)\|_\infty\leq R<R'$; by continuity of $u_{R'}$, either $\tau_{R'}=T$ or $\tau_{R'}>\tau_R$.

(2) By (1), $u$ is well defined on $[0,\tau)\times D$, and for every $t<\tau$ there is $R$ with $t<\tau_R$, so that $u=u_R$ on $[0,t]\times D$ is continuous and satisfies \eqref{eq:mild} at all $(s,x)$ with $s\leq t$, since $f_R(u_R)=f(u_R)$ there. Adaptedness follows from that of the $u_R$ and from $\{t<\tau_R\}\in\mathcal F_t$.

(3) Let $\tilde\tau_R=\inf\{t<\theta:\|\tilde u(t)\|_\infty\geq R\}\wedge\theta$. On $[0,\tilde\tau_R\wedge\tau_R)$ both $\tilde u$ and $u_R$ take values in $[-R,R]$ and satisfy \eqref{eq:mild} with $f$ replaced by $f_R$. The argument of (1) gives $\tilde u=u_R=u$ on $[0,\tilde\tau_R\wedge\tau_R)$. Since $\tilde u$ is continuous on $[0,t]$ for $t<\theta$, $\tilde\tau_R\to\theta$ as $R\to\infty$, and the claim follows by letting $R\to\infty$.

(4) Assume that $\tau<T$ and $\sup_{t<\tau}\|u(t)\|_\infty\leq M<\infty$, and let $R>M$. By (1), $\tau_R<\tau_{R+1}\leq\tau$ since $\tau_R<T$. By continuity of $u_R$ and the definition of $\tau_R$, $\|u_R(\tau_R)\|_\infty=R$ (recall $\tau_R>0$), and since $\tau_R<\tau$, $u(\tau_R)=u_R(\tau_R)$. This contradicts $\|u(\tau_R)\|_\infty\leq M<R$. For the equivalent formulation, if $\sup_{t<\tau}\|u(t)\|_\infty\leq M$ then $\tau=T$, and the same argument shows that $\tau_R<T$ is impossible for $R>M$; hence $\tau_R=T$, $u=u_R$ on $[0,T]\times D$, and $u$ satisfies \eqref{eq:mild} on $[0,T]$.
\end{proof}

\begin{remark}
Theorem~\ref{thm:local} reduces global existence to an a priori bound on $\sup_{t<\tau}\|u(t)\|_\infty$. We obtain it in Section~\ref{sec:global} by constructing directly a global solution, bounded in $L^\infty$, as the limit of spectral Galerkin approximations, and by using the uniqueness statement (3).
\end{remark}
\section{Global existence}\label{sec:global}

As in Section~\ref{sec:local}, we fix $\omega$ in the event of probability one on which Lemmas~\ref{lem:Zholder} and~\ref{lem:Zx} hold. In particular $Z\in C([0,T]\times D)$, $\partial_xZ\in C([0,T]\times D)$, $\partial_xZ=0$ on $[0,T]\times\{0,\pi\}$, and $K_Z<\infty$. The strategy is the following. We look for $u=v+Z$, where $v$ solves the random equation
\begin{equation}\label{eq:v}
\partial_t v=\Delta\big[(1+\Delta)^2v+f(v+Z)\big],\qquad v(0)=u_0,
\end{equation}
which contains no stochastic integral. We approximate \eqref{eq:v} by spectral Galerkin systems, for which the energy computations are rigorous, prove bounds uniform in the dimension, pass to the limit, and conclude with the uniqueness part of Theorem~\ref{thm:local}.

\subsection{Galerkin approximation}

For $N\geq1$ we consider the system of ordinary differential equations
\begin{equation}\label{eq:galerkin}
\xi_k'(t)=-\lambda_k\xi_k(t)-k^2\big\langle f\big(v^N(t)+Z(t)\big),e_k\big\rangle,\qquad \xi_k(0)=\hat u_{0,k},\qquad 0\leq k\leq N,
\end{equation}
where $v^N(t)=\sum_{k=0}^N\xi_k(t)e_k$. This is the projection of \eqref{eq:v} on $\mathrm{span}(e_0,\dots,e_N)$, written in weak form: $\langle\Delta\varphi,e_k\rangle$ is replaced by $\langle\varphi,\Delta e_k\rangle=-k^2\langle\varphi,e_k\rangle$, which only requires $\varphi\in L^2(D)$.

\begin{lemma}\label{lem:galerkin}
Assume that every solution of \eqref{eq:galerkin} on a subinterval $[0,t_1)\subset[0,T]$ satisfies $\sup_{t<t_1}|\xi(t)|<\infty$ (this is proved in Proposition~\ref{prop:L2}). Then \eqref{eq:galerkin} has a unique solution $\xi\in C^1([0,T];\R^{N+1})$. The function $v^N$ is continuous on $[0,T]\times D$ and adapted, $\xi_0(t)=\hat u_{0,0}$ for all $t$, and
\begin{equation}\label{eq:galerkinmild}
v^N(t,x)=S(t)P_Nu_0(x)+N^N_{f(v^N+Z)}(t,x),\qquad (t,x)\in[0,T]\times D .
\end{equation}
\end{lemma}

\begin{proof}
The right-hand side $F(t,\xi)$ of \eqref{eq:galerkin} is continuous on $[0,T]\times\R^{N+1}$, since $Z$ is continuous and $f$ is a polynomial, and it is locally Lipschitz in $\xi$, uniformly in $t\in[0,T]$: indeed, for $|\xi|,|\eta|\leq r$, $\|v_\xi\|_\infty,\|v_\eta\|_\infty\leq C_Nr$ and
\[
\big|\langle f(v_\xi+Z(t))-f(v_\eta+Z(t)),e_k\big\rangle\big|\leq C\sup_{|a|\leq C_Nr+K_Z}|f'(a)|\;\|v_\xi-v_\eta\|\leq C(N,r,K_Z)|\xi-\eta| .
\]
By the Cauchy-Lipschitz theorem, there is a unique maximal solution on $[0,t_1)$, and under the assumption of the lemma $t_1=T$ and the solution extends to $[0,T]$. (The assumption is verified in Proposition~\ref{prop:L2} below, whose proof only uses the equations \eqref{eq:galerkin} on $[0,t_1)$ and not the present lemma, so that there is no circularity.) The equation for $k=0$ is $\xi_0'=0$. For adaptedness, fix $t\in(0,T]$ and let $E_t$ be the space of functions $\zeta\in C([0,t]\times D)$ such that $\partial_x\zeta$ exists and belongs to $C([0,t]\times D)$, with the norm $\|\zeta\|_{E_t}=\sup_{[0,t]\times D}|\zeta|+\sup_{[0,t]\times D}|\partial_x\zeta|$. The map $\zeta\mapsto(\zeta,\partial_x\zeta)$ is an isometry from $E_t$ onto a closed subspace of $C([0,t]\times D)^2$, so $E_t$ is a separable Banach space. For $\zeta\in E_t$, let $\xi^\zeta$ be the solution of \eqref{eq:galerkin} on $[0,t]$ with $Z$ replaced by $\zeta$, and $\Gamma_t(\zeta)=\xi^\zeta(t)$. We check three facts.
\begin{enumerate}
\item[(i)] \emph{$\Gamma_t$ is well defined on $E_t$ and depends only on $\zeta$ on $[0,t]$.} The right-hand side of \eqref{eq:galerkin} at time $s$ involves $\zeta(s)$ only, so the solution on $[0,t]$ is determined by the restriction of $\zeta$ to $[0,t]\times D$. Its existence on the whole interval $[0,t]$ follows from the a priori bound of Proposition~\ref{prop:L2}, whose proof uses $\zeta$ only through $K_\zeta:=\|\zeta\|_{E_t}$ and $|\bar u_0|$: it gives $|\xi^\zeta(s)|^2\leq C_0(K_\zeta)$ for $s\leq t$, with $C_0$ nondecreasing in $K_\zeta$.
\item[(ii)] \emph{$\Gamma_t$ is continuous on $E_t$.} Let $\zeta,\zeta'\in E_t$ with $\|\zeta\|_{E_t},\|\zeta'\|_{E_t}\leq M$, and $\rho=C_0(M)^{1/2}$, so that $|\xi^\zeta(s)|,|\xi^{\zeta'}(s)|\leq\rho$ and $\|v^\zeta(s)\|_\infty,\|v^{\zeta'}(s)\|_\infty\leq C_N\rho$ for $s\leq t$. With $L=\sup\{|f'(a)|:|a|\leq C_N\rho+M\}$ and $|e_k|\leq(2/\pi)^{1/2}$, for every $k\leq N$,
\begin{align*}
\big|\langle f(v^\zeta+\zeta)-f(v^{\zeta'}+\zeta'),e_k\rangle\big|&\leq\pi^{1/2}\sqrt2\,L\big(\|v^\zeta-v^{\zeta'}\|_\infty+\|\zeta-\zeta'\|_\infty\big)\\
&\leq C(N,M)\big(|\xi^\zeta-\xi^{\zeta'}|+\|\zeta-\zeta'\|_{E_t}\big),
\end{align*}
so that, by \eqref{eq:galerkin} and $\lambda_k\leq N^6$, $|\xi^\zeta(s)-\xi^{\zeta'}(s)|\leq C(N,M)\int_0^s\big(|\xi^\zeta-\xi^{\zeta'}|+\|\zeta-\zeta'\|_{E_t}\big)dr$. The Gronwall lemma gives $|\Gamma_t(\zeta)-\Gamma_t(\zeta')|\leq C(N,M)\,t\,e^{C(N,M)t}\|\zeta-\zeta'\|_{E_t}$, so $\Gamma_t$ is locally Lipschitz, hence continuous and Borel measurable.
\item[(iii)] \emph{The input is $\mathcal F_t$-measurable.} By Lemma~\ref{lem:Zx}, $Z|_{[0,t]\times D}\in E_t$ almost surely. The Borel $\sigma$-field of the separable Banach space $E_t$ is generated by the point evaluations $\zeta\mapsto\zeta(s,x)$ and $\zeta\mapsto\partial_x\zeta(s,x)$, $(s,x)$ in a countable dense subset of $[0,t]\times D$, and $Z(s,x)$, $\partial_xZ(s,x)$ are $\mathcal F_s$-measurable, hence $\mathcal F_t$-measurable, for $s\leq t$. Therefore $Z|_{[0,t]\times D}$ is an $\mathcal F_t$-measurable $E_t$-valued random variable.
\end{enumerate}
By (i)-(iii), $\xi(t)=\Gamma_t(Z|_{[0,t]\times D})$ is $\mathcal F_t$-measurable. Finally, the variation of constants formula gives
\[
\xi_k(t)=e^{-\lambda_kt}\hat u_{0,k}-k^2\int_0^te^{-\lambda_k(t-s)}\big\langle f(v^N(s)+Z(s)),e_k\big\rangle\,ds,
\]
and multiplying by $e_k(x)$ and summing over $k\leq N$ gives \eqref{eq:galerkinmild}, since
\[
\Delta_yG_N(r,x,y)=-\sum_{k\leq N}k^2e^{-\lambda_kr}e_k(x)e_k(y).
\]
\end{proof}

\subsection{An interpolation inequality with the neutral mode}

\begin{lemma}\label{lem:interp}
For every $\delta>0$ there exists $C_\delta$ such that, for every $w=\sum_{k\leq N}\hat w_ke_k$, $N\geq1$,
\[
\|\partial_xw\|^2\leq\delta\,\|A^{1/2}w\|^2+C_\delta\,\|w\|^2 .
\]
\end{lemma}

\begin{proof}
By \eqref{eq:basis}, $\partial_xw=-\sum_{k\geq1}k\,\hat w_kg_k$, so $\|\partial_xw\|^2=\sum_{k\geq1}k^2|\hat w_k|^2$. Let $M\geq4$ with $\frac{16}{9M^2}\leq\delta$. The modes with $k^2\leq M$, among which the neutral mode $k=1$, contribute at most $M\|w\|^2$. For $k^2>M$ we have $k\geq2$, and Lemma~\ref{lem:spectrum} gives
\[
k^2=\frac{k^6}{k^4}\leq\frac{16}9\,\frac{\lambda_k}{k^4}\leq\frac{16}{9M^2}\,\lambda_k\leq\delta\lambda_k .
\]
Summing, $\|\partial_xw\|^2\leq M\|w\|^2+\delta\|A^{1/2}w\|^2$.
\end{proof}

\begin{remark}
Lemma~\ref{lem:interp} is the first of the two places where the neutral mode $e_1$ requires an argument specific to the PFC equation; the second is the $H^{-1}$ estimate of Lemma~\ref{lem:Hminus1}, where the neutral mode is controlled by the quartic term. Since $\lambda_1=0$, the dissipation $\|A^{1/2}w\|^2$ does not control $\hat w_1$, and the corresponding part of $\|\partial_xw\|^2$ is controlled by $\|w\|^2$ instead, whose time integral is controlled by the $H^{-1}$ estimate of Lemma~\ref{lem:Hminus1}, where the neutral mode is in turn controlled by the quartic term. For Cahn-Hilliard, $\|\partial_x\Delta w\|^2$ controls $\|\partial_xw\|^2$ directly.
\end{remark}

\subsection{An \texorpdfstring{$H^{-1}$}{H-1} estimate exploiting the quartic potential}

Let $\bar u_0=\pi^{-1}\int_Du_0\,dx$ be the mean of $u_0$, and, for $w\in L^2(D)$, $\|w\|^2_{-1}:=\sum_{k\geq1}k^{-2}|\hat w_k|^2$; for $w$ with zero mean this is the norm of $H^{-1}$. We also set
\[
m_0:=K_Z+|\bar u_0| .
\]
The following estimate uses the growth of the quartic potential rather than the linear dissipation. It is the analogue, for our conservative equation, of the $H^{-1}$ estimate of Da Prato and Debussche \cite{DPD96} for the Cahn-Hilliard equation, and of the $L^p$ estimates of Doering used by Bertini, Brassesco, Butt\`a and Presutti \cite{BBBP02}, which also exploit the dissipation produced by the quartic term. Its point is that it gives bounds which are polynomial in $K_Z$.

\begin{lemma}\label{lem:Hminus1}
Let $\xi$ be a solution of \eqref{eq:galerkin} on $[0,t_1)\subset[0,T]$ and $w^N=v^N+Z$. Then, for $t<t_1$,
\begin{equation}\label{eq:Hm1}
\frac12\|v^N(t)-\bar u_0\|^2_{-1}+\int_0^t\sum_{k=1}^N(1-k^2)^2\xi_k(s)^2\,ds+\frac14\int_0^t\|w^N(s)\|^4_{L^4}\,ds\leq\frac12\|u_0\|^2+C_\varepsilon T\big(1+m_0^4\big),
\end{equation}
with $C_\varepsilon=\pi\big(\varepsilon^2+7+|\varepsilon|^{4/3}\big)$. Consequently
\begin{equation}\label{eq:L2int}
\int_0^t\|v^N(s)\|^2\,ds\leq C_7:=\frac{\pi T}2+8\|u_0\|^2+16C_\varepsilon T\big(1+m_0^4\big)+4\pi TK_Z^4 .
\end{equation}
\end{lemma}

\begin{proof}
Throughout the proof we write $v=v^N$, $w=w^N$ and $\varphi=f(w)$.

\emph{Step 1: identity.} By Lemma~\ref{lem:galerkin}, $\xi_0$ is constant, and since $\xi_0(0)=\hat u_{0,0}=\pi^{-1/2}\int_Du_0$ and $e_0=\pi^{-1/2}$, we have $\xi_0e_0=\bar u_0$; hence $v-\bar u_0=\sum_{k=1}^N\xi_ke_k$ and $\|v-\bar u_0\|^2_{-1}=\sum_{k=1}^N\xi_k^2/k^2$. By \eqref{eq:galerkin}, since $\lambda_k/k^2=(1-k^2)^2$,
\begin{align*}
\frac12\frac{d}{dt}\|v-\bar u_0\|^2_{-1}
&=\sum_{k=1}^N\frac{\xi_k\,\xi_k'}{k^2}
=-\sum_{k=1}^N\frac{\lambda_k}{k^2}\,\xi_k^2-\sum_{k=1}^N\xi_k\,\langle\varphi,e_k\rangle\\
&=-\sum_{k=1}^N(1-k^2)^2\xi_k^2-\Big\langle\varphi,\sum_{k=1}^N\xi_ke_k\Big\rangle
=-\sum_{k=1}^N(1-k^2)^2\xi_k^2-\langle\varphi,v-\bar u_0\rangle .
\end{align*}
Note that the neutral mode $k=1$ does not contribute to the first sum; it will be controlled by the quartic term.

\emph{Step 2: lower bound on the nonlinear term.} Let $\zeta=Z+\bar u_0$, so that $v-\bar u_0=w-\zeta$ and $|\zeta|\leq m_0$ on $[0,T]\times D$. Then
\[
\langle\varphi,v-\bar u_0\rangle=\int_D\big(w^3-\varepsilon w\big)(w-\zeta)\,dx=\int_Dw^4\,dx-\varepsilon\!\int_Dw^2\,dx-\int_Dw^3\zeta\,dx+\varepsilon\!\int_Dw\,\zeta\,dx .
\]
We bound the last three terms pointwise. By the inequality $|\varepsilon|a^2\leq\frac14a^4+\varepsilon^2$, by Young's inequality with exponents $\frac43$ and $4$ in the forms $|a|^3m\leq\frac14a^4+\frac{27}4m^4$ and $|a|\,b\leq\frac14a^4+\frac34b^{4/3}$,
\begin{align*}
\varepsilon w^2&\leq\tfrac14w^4+\varepsilon^2,\\
|w^3\zeta|&\leq|w|^3m_0\leq\tfrac14w^4+\tfrac{27}4m_0^4,\\
|\varepsilon w\zeta|&\leq|w|\,|\varepsilon|m_0\leq\tfrac14w^4+\tfrac34|\varepsilon|^{4/3}m_0^{4/3}\leq\tfrac14w^4+|\varepsilon|^{4/3}\big(1+m_0^4\big).
\end{align*}
Integrating over $D$, with $|D|=\pi$,
\begin{align*}
\langle\varphi,v-\bar u_0\rangle
&\geq\int_Dw^4-\frac34\int_Dw^4-\pi\Big(\varepsilon^2+\frac{27}4m_0^4+|\varepsilon|^{4/3}\big(1+m_0^4\big)\Big)\\
&\geq\frac14\|w\|^4_{L^4}-C_\varepsilon\big(1+m_0^4\big),
\end{align*}
with $C_\varepsilon=\pi(\varepsilon^2+7+|\varepsilon|^{4/3})$.

\emph{Step 3: integration.} Combining Steps 1 and 2,
\[
\frac12\frac{d}{dt}\|v-\bar u_0\|^2_{-1}+\sum_{k=1}^N(1-k^2)^2\xi_k^2+\frac14\|w\|^4_{L^4}\leq C_\varepsilon\big(1+m_0^4\big),
\]
and integrating on $[0,t]$ with $\|v(0)-\bar u_0\|^2_{-1}=\sum_{k=1}^N\hat u_{0,k}^2/k^2\leq\|u_0\|^2$ gives \eqref{eq:Hm1}.

\emph{Step 4: proof of \eqref{eq:L2int}.} By the Cauchy-Schwarz inequality and $ab\leq\frac12a^2+\frac12b^2$,
\[
\|v\|^2\leq|D|^{1/2}\|v\|^2_{L^4}\leq\frac\pi2+\frac12\|v\|^4_{L^4},
\]
and, since $|a+b|^4\leq8(a^4+b^4)$ and $v=w-Z$, $\|v\|^4_{L^4}\leq8\|w\|^4_{L^4}+8\pi K_Z^4$. Hence, by \eqref{eq:Hm1},
\begin{align*}
\int_0^t\|v(s)\|^2\,ds
&\leq\frac{\pi T}2+4\int_0^t\|w(s)\|^4_{L^4}\,ds+4\pi TK_Z^4\\
&\leq\frac{\pi T}2+16\Big(\frac12\|u_0\|^2+C_\varepsilon T\big(1+m_0^4\big)\Big)+4\pi TK_Z^4=C_7 .
\end{align*}
\end{proof}

\subsection{Uniform \texorpdfstring{$L^2$}{L2} estimate}

\begin{proposition}\label{prop:L2}
Let $\xi$ be a solution of \eqref{eq:galerkin} on $[0,t_1)\subset[0,T]$. Then, for $t<t_1$,
\begin{equation}\label{eq:L2}
\|v^N(t)\|^2+\int_0^t\|A^{1/2}v^N(s)\|^2\,ds\leq C_0:=\|u_0\|^2+C_1\big(1+K_Z^2\big)C_7+C_2T\big(1+K_Z^4\big),
\end{equation}
where $C_7$ is given by \eqref{eq:L2int} and $C_1,C_2$ depend only on $\varepsilon$. In particular $C_0$ does not depend on $N$ nor on $t_1$, and it is bounded by a polynomial in $\|u_0\|$ and $m_0=K_Z+|\bar u_0|$, with coefficients depending only on $T$ and $\varepsilon$.
\end{proposition}

\begin{proof}
Throughout the proof we write $v=v^N$ and $\varphi=f(v+Z)$, and $t<t_1$.

\emph{Step 1: energy identity.} Since $\xi\in C^1$ and $\|v(t)\|^2=\sum_{k\leq N}\xi_k(t)^2$, \eqref{eq:galerkin} gives
\begin{align*}
\frac12\frac{d}{dt}\|v\|^2
&=\sum_{k\leq N}\xi_k\,\xi_k'
=-\sum_{k\leq N}\lambda_k\xi_k^2-\sum_{k\leq N}k^2\xi_k\,\langle\varphi,e_k\rangle\\
&=-\|A^{1/2}v\|^2+\Big\langle\varphi,-\sum_{k\leq N}k^2\xi_ke_k\Big\rangle
=-\|A^{1/2}v\|^2+\langle\varphi,\Delta v\rangle .
\end{align*}
Now $v(t,\cdot)$ is a trigonometric polynomial in the $e_k$, so $\partial_xv(t,0)=\partial_xv(t,\pi)=0$, and $\varphi(t,\cdot)\in C^1(D)$ by Lemma~\ref{lem:Zx}, with $\partial_x\varphi=f'(v+Z)(\partial_xv+\partial_xZ)$ and $f'(a)=3a^2-\varepsilon$. An integration by parts, whose boundary term $[\varphi\,\partial_xv]_0^\pi$ vanishes, gives
\begin{align*}
\langle\varphi,\Delta v\rangle
&=-\langle\partial_x\varphi,\partial_xv\rangle
=-\int_D\big(3(v+Z)^2-\varepsilon\big)\big(\partial_xv+\partial_xZ\big)\partial_xv\,dx\\
&=-3\!\int_D(v+Z)^2|\partial_xv|^2\,dx-3\!\int_D(v+Z)^2\,\partial_xZ\,\partial_xv\,dx
+\varepsilon\|\partial_xv\|^2+\varepsilon\langle\partial_xZ,\partial_xv\rangle .
\end{align*}

\emph{Step 2: the cubic part has a good sign.} By Young's inequality $|ab|\leq\frac12a^2+\frac12b^2$ applied pointwise with $a=\partial_xv$, $b=\partial_xZ$,
\begin{align*}
-3\!\int_D(v+Z)^2|\partial_xv|^2-3\!\int_D(v+Z)^2\,\partial_xZ\,\partial_xv
&\leq-3\!\int_D(v+Z)^2|\partial_xv|^2\\
&\qquad+\frac32\!\int_D(v+Z)^2|\partial_xv|^2+\frac32\!\int_D(v+Z)^2|\partial_xZ|^2\\
&=-\frac32\!\int_D(v+Z)^2|\partial_xv|^2+\frac32\!\int_D(v+Z)^2|\partial_xZ|^2\\
&\leq\frac32\,\|\partial_xZ(t)\|_\infty^2\int_D(v+Z)^2\,dx\\
&\leq\frac32K_Z^2\int_D2\big(v^2+Z^2\big)\,dx\\
&\leq3K_Z^2\|v\|^2+3\pi K_Z^4 .
\end{align*}
Only $\partial_xZ$, and not $\partial_x^2Z$, appears in this bound.

\emph{Step 3: the linear part.} By the Cauchy-Schwarz and Young inequalities,
\[
\varepsilon\|\partial_xv\|^2+\varepsilon\langle\partial_xZ,\partial_xv\rangle\leq|\varepsilon|\|\partial_xv\|^2+\frac12\|\partial_xv\|^2+\frac{\varepsilon^2}2\|\partial_xZ\|^2\leq\Big(|\varepsilon|+\frac12\Big)\|\partial_xv\|^2+\frac{\pi\varepsilon^2}2K_Z^2 .
\]
By Lemma~\ref{lem:interp} with $\delta=(2|\varepsilon|+1)^{-1}$, so that $(|\varepsilon|+\frac12)\delta=\frac12$,
\[
\Big(|\varepsilon|+\frac12\Big)\|\partial_xv\|^2\leq\frac12\|A^{1/2}v\|^2+C_\varepsilon\|v\|^2,\qquad C_\varepsilon:=\Big(|\varepsilon|+\frac12\Big)C_\delta .
\]

\emph{Step 4: differential inequality.} Collecting Steps 1 to 3,
\begin{align*}
\frac12\frac{d}{dt}\|v\|^2
&\leq-\|A^{1/2}v\|^2+3K_Z^2\|v\|^2+3\pi K_Z^4+\frac12\|A^{1/2}v\|^2+C_\varepsilon\|v\|^2+\frac{\pi\varepsilon^2}2K_Z^2\\
&=-\frac12\|A^{1/2}v\|^2+\big(3K_Z^2+C_\varepsilon\big)\|v\|^2+3\pi K_Z^4+\frac{\pi\varepsilon^2}2K_Z^2 .
\end{align*}
Multiplying by $2$ and using $K_Z^2\leq1+K_Z^4$, we obtain
\begin{equation}\label{eq:diffineq}
\frac{d}{dt}\|v\|^2+\|A^{1/2}v\|^2\leq C_1\big(1+K_Z^2\big)\|v\|^2+C_2\big(1+K_Z^4\big),
\end{equation}
with $C_1=\max(6,2C_\varepsilon)$ and $C_2=6\pi+\pi\varepsilon^2$, which depend only on $\varepsilon$.

\emph{Step 5: integration without Gronwall.} Instead of applying the Gronwall lemma to \eqref{eq:diffineq}, which would produce a factor $e^{C_1T(1+K_Z^2)}$, we integrate it on $[0,t]$ and use the bound \eqref{eq:L2int} of Lemma~\ref{lem:Hminus1} on the time integral of $\|v\|^2$. Since $\|v(0)\|=\|P_Nu_0\|\leq\|u_0\|$,
\begin{align*}
\|v(t)\|^2+\int_0^t\|A^{1/2}v(s)\|^2\,ds
&\leq\|v(0)\|^2+C_1\big(1+K_Z^2\big)\int_0^t\|v(s)\|^2\,ds+C_2\,t\,\big(1+K_Z^4\big)\\
&\leq\|u_0\|^2+C_1\big(1+K_Z^2\big)\,C_7+C_2T\big(1+K_Z^4\big)=C_0,
\end{align*}
which is \eqref{eq:L2}. The polynomial character of $C_0$ follows from that of $C_7$.
\end{proof}

By Proposition~\ref{prop:L2}, $|\xi(t)|=\|v^N(t)\|\leq C_0^{1/2}$ on any interval of existence, so the assumption of Lemma~\ref{lem:galerkin} holds: for every $N$, $v^N$ is defined on $[0,T]$ and satisfies \eqref{eq:L2} with $t\leq T$.

\subsection{Uniform \texorpdfstring{$L^\infty$}{L-infinity} estimate}

\begin{proposition}\label{prop:Linf}
There exists $C'_0<\infty$, depending only on $T$, $\varepsilon$, $\|u_0\|_{H^1}$ and $K_Z$, such that $\sup_{N\geq1}\sup_{t\in[0,T]}\|v^N(t)\|_\infty\leq C'_0$. Here $C'_0\leq P(\|u_0\|_{H^1},K_Z)$ for a polynomial $P$ whose coefficients depend only on $T$ and $\varepsilon$.
\end{proposition}

\begin{proof}
\emph{Step 1: an $H^3$ bound.} For $k\geq2$, $1+k^2\leq\frac54k^2$, hence $(1+k^2)^3\leq\frac{125}{64}k^6\leq\frac{125}{64}\cdot\frac{16}9\lambda_k=\frac{125}{36}\lambda_k\leq4\lambda_k$ by Lemma~\ref{lem:spectrum}; for $k\in\{0,1\}$, $(1+k^2)^3\leq8$. Therefore, for every $w=\sum_k\hat w_ke_k$,
\[
\|w\|^2_{H^3}=\sum_k(1+k^2)^3|\hat w_k|^2\leq8\big(|\hat w_0|^2+|\hat w_1|^2\big)+4\sum_{k\geq2}\lambda_k|\hat w_k|^2\leq8\|w\|^2+4\|A^{1/2}w\|^2,
\]
and \eqref{eq:L2} gives
\begin{equation}\label{eq:H3}
\int_0^T\|v^N(s)\|^2_{H^3}\,ds\leq8\int_0^T\|v^N(s)\|^2ds+4\int_0^T\|A^{1/2}v^N(s)\|^2ds\leq(8T+4)\,C_0 .
\end{equation}

\emph{Step 2: interpolation.} Let $w=\sum_k\hat w_ke_k$. By \eqref{eq:basis} and the Cauchy-Schwarz inequality,
\begin{align*}
\|w\|_\infty&\leq\Big(\frac2\pi\Big)^{1/2}\sum_k|\hat w_k|
=\Big(\frac2\pi\Big)^{1/2}\sum_k(1+k^2)^{1/3}|\hat w_k|\,(1+k^2)^{-1/3}\\
&\leq\Big(\frac2\pi\Big)^{1/2}\Big(\sum_k(1+k^2)^{2/3}|\hat w_k|^2\Big)^{1/2}\Big(\sum_k(1+k^2)^{-2/3}\Big)^{1/2}=:c_1\,\|w\|_{H^{2/3}},
\end{align*}
where $c_1<\infty$ because $\frac43>1$. Next, writing $(1+k^2)^{2/3}|\hat w_k|^2=|\hat w_k|^{14/9}\big((1+k^2)^3|\hat w_k|^2\big)^{2/9}$ and applying the H\"older inequality with exponents $\frac97$ and $\frac92$,
\[
\|w\|^2_{H^{2/3}}\leq\Big(\sum_k|\hat w_k|^2\Big)^{7/9}\Big(\sum_k(1+k^2)^3|\hat w_k|^2\Big)^{2/9}=\|w\|^{14/9}\,\|w\|_{H^3}^{4/9}.
\]
Hence $\|w\|_\infty\leq c_1\|w\|^{7/9}\|w\|^{2/9}_{H^3}$, and raising to the power $9$,
\[
\|v^N(s)\|^9_\infty\leq c_1^9\,\|v^N(s)\|^{7}\,\|v^N(s)\|^{2}_{H^3}\leq c_1^9\,C_0^{7/2}\,\|v^N(s)\|^{2}_{H^3}.
\]
Integrating in $s$ and using \eqref{eq:H3},
\begin{equation}\label{eq:L9}
\int_0^T\|v^N(s)\|^9_\infty\,ds\leq c_1^9\,C_0^{7/2}(8T+4)\,C_0=:C_3 .
\end{equation}

\emph{Step 3: mild formulation.} Let $\varphi^N=f(v^N+Z)$. Since $|f(a)|\leq|a|^3+|\varepsilon||a|\leq2|a|^3+|\varepsilon|^{3/2}$ (by Young's inequality) and $|a+b|^3\leq4(|a|^3+|b|^3)$,
\[
\|\varphi^N(s)\|_\infty\leq8\big(\|v^N(s)\|^3_\infty+K_Z^3\big)+|\varepsilon|^{3/2}=:g_N(s),
\]
and, by $(a+b+c)^3\leq9(a^3+b^3+c^3)$ and \eqref{eq:L9},
\begin{align*}
\int_0^Tg_N(s)^3\,ds&\leq9\Big(8^3\int_0^T\|v^N(s)\|^9_\infty ds+8^3TK_Z^9+T|\varepsilon|^{9/2}\Big)\\
&\leq9\big(512\,C_3+512\,TK_Z^9+T|\varepsilon|^{9/2}\big)=:C_4 .
\end{align*}
By \eqref{eq:galerkinmild}, \eqref{eq:initial}, \eqref{eq:conv1} for $N^N$, and the H\"older inequality with exponents $\frac32$ and $3$,
\begin{align*}
\|v^N(t)\|_\infty
&\leq\|S(t)P_Nu_0\|_\infty+\|N^N_{\varphi^N}(t)\|_\infty
\leq C\|u_0\|_{H^1}+C\int_0^t\big(1+(t-s)^{-1/2}\big)g_N(s)\,ds\\
&\leq C\|u_0\|_{H^1}+C\Big(\int_0^t\big(1+(t-s)^{-1/2}\big)^{3/2}ds\Big)^{2/3}\Big(\int_0^tg_N(s)^3\,ds\Big)^{1/3}\\
&\leq C\|u_0\|_{H^1}+C\Big(\sqrt2\int_0^T\big(1+r^{-3/4}\big)dr\Big)^{2/3}C_4^{1/3}\\
&=C\|u_0\|_{H^1}+C\big(\sqrt2\,(T+4T^{1/4})\big)^{2/3}C_4^{1/3}=:C'_0,
\end{align*}
where we used $(a+b)^{3/2}\leq\sqrt2\,(a^{3/2}+b^{3/2})$. The constant $C'_0$ depends only on $T$, $\varepsilon$, $\|u_0\|_{H^1}$ and $K_Z$, and not on $N$ or $t$. Since $C_0$ is polynomial in $\|u_0\|$ and $m_0$ (Proposition~\ref{prop:L2}), and $|\bar u_0|\leq\pi^{-1/2}\|u_0\|$, the constants $C_3$, $C_4$ and $C'_0$ are polynomial in $\|u_0\|_{H^1}$ and $K_Z$.
\end{proof}

\subsection{Passage to the limit}

\begin{proposition}\label{prop:limit}
The sequence $(v^N)$ converges uniformly on $[0,T]\times D$ to a function $\tilde v\in C([0,T]\times D)$, which is adapted, satisfies $\sup_{t}\|\tilde v(t)\|_\infty\leq C'_0$, and
\begin{equation}\label{eq:vtilde}
\tilde v(t,x)=S(t)u_0(x)+N_{f(\tilde v+Z)}(t,x),\qquad (t,x)\in[0,T]\times D .
\end{equation}
\end{proposition}

\begin{proof}
\emph{Step 1: bounds on the nonlinearity.} Let $M=C'_0+K_Z$, $L_M=\sup_{|a|\leq M}|f'(a)|=3M^2+|\varepsilon|$ and $B_M=\sup_{|a|\leq M}|f(a)|\leq M^3+|\varepsilon|M$. By Proposition~\ref{prop:Linf}, $\|v^N(t)+Z(t)\|_\infty\leq M$ for all $N$ and $t$, so that $\varphi^N=f(v^N+Z)$ satisfies
\[
\|\varphi^N\|_\infty\leq B_M,\qquad|\varphi^N(s,y)-\varphi^{N'}(s,y)|\leq L_M\,|v^N(s,y)-v^{N'}(s,y)| .
\]

\emph{Step 2: Cauchy estimate.} Let $N<N'$. By \eqref{eq:galerkinmild},
\begin{align*}
v^N-v^{N'}
&=S(\cdot)P_Nu_0-S(\cdot)P_{N'}u_0+N^N_{\varphi^N}-N^{N'}_{\varphi^{N'}}\\
&=S(\cdot)(P_N-P_{N'})u_0+\big(N^N_{\varphi^N}-N_{\varphi^N}\big)-\big(N^{N'}_{\varphi^N}-N_{\varphi^N}\big)+N^{N'}_{\varphi^N-\varphi^{N'}} .
\end{align*}
We bound the four terms. By the proof of Lemma~\ref{lem:initial}, the first is bounded by $\eta_N:=C\big(\sum_{k>N}(1+k^2)|\hat u_{0,k}|^2\big)^{1/2}$. By \eqref{eq:conv4}, each of the second and third terms is bounded by $\eta'_N:=2B_M\sum_{k>N}k^2/\lambda_k$. By \eqref{eq:conv1} applied to $N^{N'}$ and Step~1, the fourth term at time $t$ is bounded by $CL_M\int_0^t(1+(t-s)^{-1/2})\|v^N(s)-v^{N'}(s)\|_\infty\,ds$. Hence, for all $t\in[0,T]$,
\[
\|v^N(t)-v^{N'}(t)\|_\infty\leq\eta_N+2\eta'_N+CL_M\int_0^t\big(1+(t-s)^{-1/2}\big)\|v^N(s)-v^{N'}(s)\|_\infty\,ds,
\]
and Lemma~\ref{lem:gronwall} with $a=\eta_N+2\eta'_N$ and $b=CL_M$ gives
\[
\sup_{t\in[0,T]}\|v^N(t)-v^{N'}(t)\|_\infty\leq2e^{\gamma(CL_M)T}\big(\eta_N+2\eta'_N\big)\xrightarrow[N\to\infty]{}0,
\]
since $\eta_N\to0$ (tail of the convergent series $\sum(1+k^2)|\hat u_{0,k}|^2$) and $\eta'_N\to0$ (tail of $\sum k^2/\lambda_k$).

\emph{Step 3: limit.} Hence $(v^N)$ is a Cauchy sequence in $C([0,T]\times D)$, which converges uniformly to some $\tilde v\in C([0,T]\times D)$ with $\|\tilde v\|_\infty\leq C'_0$. Each $v^N(t,x)$ is $\mathcal F_t$-measurable by Lemma~\ref{lem:galerkin}, and so is the pointwise limit $\tilde v(t,x)$.

\emph{Step 4: equation for the limit.} Let $\tilde\varphi=f(\tilde v+Z)$, so that $\|\varphi^N-\tilde\varphi\|_\infty\leq L_M\|v^N-\tilde v\|_\infty$. By \eqref{eq:galerkinmild},
\begin{align*}
\big|\tilde v-S(\cdot)u_0-N_{\tilde\varphi}\big|
&\leq|\tilde v-v^N|+|S(\cdot)(P_N-I)u_0|+|N^N_{\varphi^N}-N_{\varphi^N}|+|N_{\varphi^N-\tilde\varphi}|\\
&\leq\|\tilde v-v^N\|_\infty+\eta_N+\eta'_N+C\big(T+2T^{1/2}\big)L_M\|v^N-\tilde v\|_\infty ,
\end{align*}
by \eqref{eq:initial}, \eqref{eq:conv4} and \eqref{eq:conv1}. The right-hand side tends to $0$ as $N\to\infty$, which proves \eqref{eq:vtilde}.
\end{proof}

\subsection{Main result}

\begin{theorem}[Global well-posedness]\label{thm:global}
Let $u_0\in H^1(D)$, $\varepsilon,\sigma\in\R$ and $T>0$. Almost surely, \eqref{eq:pfc}-\eqref{eq:bc} has a unique mild solution $u$ on $[0,T]$ in the sense of Definition~\ref{def:mild}. It satisfies $\int_Du(t,x)\,dx=\int_Du_0(x)\,dx$ for all $t\in[0,T]$ and
\begin{equation}\label{eq:globalbound}
\sup_{(t,x)\in[0,T]\times D}|u(t,x)|\leq C'_0+K_Z<\infty .
\end{equation}
Moreover the time $\tau$ of Theorem~\ref{thm:local} equals $T$, $u=u_R$ on $[0,T]\times D$ for every $R>C'_0+K_Z$, and
\begin{equation}\label{eq:moments}
\E\Big[\sup_{(t,x)\in[0,T]\times D}|u(t,x)|^p\Big]<\infty\qquad\text{for every }p\geq1 .
\end{equation}
\end{theorem}

\begin{proof}
Let $\tilde u=\tilde v+Z$. By Proposition~\ref{prop:limit}, $\tilde u$ is continuous and adapted, and \eqref{eq:vtilde} is exactly \eqref{eq:mild} on $[0,T]$; hence $\tilde u$ is a mild solution on $[0,T]$, and it satisfies \eqref{eq:globalbound}. By Theorem~\ref{thm:local}(3), $\tilde u=u$ on $[0,\tau)$, so that $\sup_{t<\tau}\|u(t)\|_\infty\leq C'_0+K_Z$. Theorem~\ref{thm:local}(4) then gives $\tau=T$, $u=u_R$ on $[0,T]\times D$ for $R>C'_0+K_Z$, and $u=\tilde u$. Uniqueness among mild solutions on $[0,T]$ follows from Theorem~\ref{thm:local}(3), and the conservation of the mean from Proposition~\ref{prop:mass}. Finally, by Proposition~\ref{prop:Linf}, $\sup_{t,x}|u(t,x)|\leq P(\|u_0\|_{H^1},K_Z)+K_Z$ with $P$ a polynomial, and $K_Z$ has moments of all orders by Lemma~\ref{lem:Zx}; this gives \eqref{eq:moments}.
\end{proof}

\begin{remark}[Moments]\label{rem:moments}
The moment bound \eqref{eq:moments} rests on Lemma~\ref{lem:Hminus1}. A direct Gronwall argument on the $L^2$ estimate \eqref{eq:diffineq} would only give a bound of order $\exp(CK_Z^2)$, which has no finite moments of high order since $K_Z$ has Gaussian tails. The $H^{-1}$ estimate avoids this by using the quartic growth of the potential, in the spirit of \cite{DPD96,BBBP02}. The existence of a density in Section~\ref{sec:density} is proved by localization and does not use \eqref{eq:moments}; the moments are however a necessary ingredient for quantitative questions such as the smoothness of the density.
\end{remark}
\section{H\"older regularity of the solution}\label{sec:regularity}

\begin{lemma}[Regularity of the initial term]\label{lem:initialreg}
Let $u_0\in H^1(D)$. For all $s,t\in[0,T]$ and $x,y\in D$,
\[
|S(t)u_0(x)-S(t)u_0(y)|\leq C\|u_0\|_{H^1}|x-y|^{1/2},\qquad |S(t)u_0(x)-S(s)u_0(x)|\leq C\|u_0\|_{H^1}|t-s|^{1/12}.
\]
Moreover, for every $t_0>0$, $S(\cdot)u_0$ is Lipschitz in $t$ and $C^\infty$ in $x$ on $[t_0,T]\times D$, and $\partial_xS(t)u_0$ vanishes at $x\in\{0,\pi\}$ for $t>0$.
\end{lemma}

\begin{proof}
Let $a_k=(1+k^2)^{1/2}|\hat u_{0,k}|$, so that $\sum_ka_k^2=\|u_0\|^2_{H^1}$.

\emph{Space.} Let $h=|x-y|\in(0,\pi]$. Since $|e_k(x)-e_k(y)|\leq(2/\pi)^{1/2}\min(kh,2)$ by the mean value theorem, the Cauchy-Schwarz inequality and Lemma~\ref{lem:sums}(c) with $p=2$ give
\begin{align*}
|S(t)u_0(x)-S(t)u_0(y)|
&\leq\sum_{k\geq1}e^{-\lambda_kt}|\hat u_{0,k}|\,|e_k(x)-e_k(y)|
\leq\Big(\frac2\pi\Big)^{1/2}\sum_{k\geq1}a_k\,\frac{\min(kh,2)}{(1+k^2)^{1/2}}\\
&\leq\Big(\frac2\pi\Big)^{1/2}\Big(\sum_{k\geq1}a_k^2\Big)^{1/2}\Big(\sum_{k\geq1}\frac{\min(k^2h^2,4)}{1+k^2}\Big)^{1/2}\\
&\leq\Big(\frac2\pi\Big)^{1/2}\|u_0\|_{H^1}\Big(\sum_{k\geq1}k^{-2}\min(k^2h^2,4)\Big)^{1/2}
\leq\Big(\frac2\pi\Big)^{1/2}C_2^{1/2}\,\|u_0\|_{H^1}\,h^{1/2}.
\end{align*}

\emph{Time.} Let $0\leq s<t\leq T$ and $h=t-s$; for $h\geq1$ the bound follows from \eqref{eq:initial}, so let $h<1$. The neutral modes do not depend on time. For $k\geq2$, $|e^{-\lambda_kt}-e^{-\lambda_ks}|=e^{-\lambda_ks}(1-e^{-\lambda_kh})\leq\min(1,\lambda_kh)\leq\min(1,k^6h)$, using $\lambda_k\leq k^6$. Let $N=\lfloor h^{-1/6}\rfloor\geq1$, so that $N\leq h^{-1/6}$ and $N+1\geq h^{-1/6}$. Then
\begin{align*}
|S(t)u_0(x)-S(s)u_0(x)|
&\leq\Big(\frac2\pi\Big)^{1/2}\sum_{k\geq2}a_k\,\frac{\min(1,k^6h)}{(1+k^2)^{1/2}}
\leq\Big(\frac2\pi\Big)^{1/2}\|u_0\|_{H^1}\Big(\sum_{k\geq2}\frac{\min(1,k^{12}h^2)}{1+k^2}\Big)^{1/2}\\
&\leq\Big(\frac2\pi\Big)^{1/2}\|u_0\|_{H^1}\Big(h^2\sum_{2\leq k\leq N}k^{10}+\sum_{k>N}k^{-2}\Big)^{1/2}\\
&\leq\Big(\frac2\pi\Big)^{1/2}\|u_0\|_{H^1}\Big(h^2N^{11}+\frac1N\Big)^{1/2}\\
&\leq\Big(\frac2\pi\Big)^{1/2}\|u_0\|_{H^1}\Big(h^{2-\frac{11}6}+\frac2{N+1}\Big)^{1/2}\\
&\leq\Big(\frac2\pi\Big)^{1/2}\|u_0\|_{H^1}\big(h^{1/6}+2h^{1/6}\big)^{1/2}\\
&=\Big(\frac6\pi\Big)^{1/2}\|u_0\|_{H^1}\,h^{1/12},
\end{align*}
where we used $\sum_{k>N}k^{-2}\leq1/N\leq2/(N+1)$ for $N\geq1$.

\emph{Smoothness.} Let $t\geq t_0>0$ and $m,n\geq0$. The series obtained by differentiating $S(t)u_0(x)$ $m$ times in $t$ and $n$ times in $x$ has terms bounded by
\[
\lambda_k^m\,k^n\,e^{-\lambda_kt}\,|\hat u_{0,k}|\,\Big(\frac2\pi\Big)^{1/2}\leq\Big(\frac2\pi\Big)^{1/2}\|u_0\|\,k^{6m+n}e^{-\frac9{16}k^6t_0}\qquad(k\geq2),
\]
which is summable, so these series converge uniformly on $[t_0,T]\times D$. Hence $S(\cdot)u_0\in C^\infty([t_0,T]\times D)$, and in particular it is Lipschitz in $t$ and $x$ there. The derivative $\partial_xS(t)u_0=-\sum_kk\,e^{-\lambda_kt}\hat u_{0,k}g_k$ vanishes termwise at $0$ and $\pi$.
\end{proof}

\begin{lemma}[$C^1$ regularity of the deterministic convolution]\label{lem:convC1}
Let $\varphi\in L^\infty([0,T]\times D)$. For every $t\in(0,T]$, $N_\varphi(t,\cdot)\in C^1(D)$,
\[
\partial_xN_\varphi(t,x)=\int_0^t\!\!\int_D\partial_x\Delta_yG(t-s,x,y)\,\varphi(s,y)\,dy\,ds,
\]
and $\partial_xN_\varphi(t,0)=\partial_xN_\varphi(t,\pi)=0$. Also $\partial_xN_\varphi(t,\cdot)$ is Lipschitz, with $|\partial_xN_\varphi(t,x)-\partial_xN_\varphi(t,x')|\leq C\|\varphi\|_\infty|x-x'|$ for a constant $C=C(T)$.
\end{lemma}

\begin{proof}
By Lemma~\ref{lem:green} with $(i,j)=(1,2)$, $|\partial_x\Delta_yG(r,x,y)|\leq C(1+r^{-2/3})$, which is integrable on $[0,T]\times D$ uniformly in $x$. For each $r>0$, $x\mapsto\partial_x\Delta_yG(r,x,y)$ is continuous. Differentiation under the integral sign (dominated convergence applied to the difference quotients, which are bounded by $C(1+r^{-2/3})\|\varphi\|_\infty$ by the mean value theorem) gives the formula, and dominated convergence again gives the continuity of $\partial_xN_\varphi(t,\cdot)$. Finally, by \eqref{eq:basis}, $\partial_x\Delta_yG(r,x,y)=\sum_{k\geq1}k^3e^{-\lambda_kr}g_k(x)e_k(y)$ for $r>0$, which vanishes at $x\in\{0,\pi\}$. For the Lipschitz property, Lemma~\ref{lem:green} with $(i,j)=(2,2)$ gives $|\partial_x^2\Delta_yG(r,x,y)|\leq C_{2,2}(1+r^{-5/6})$, so that by the mean value theorem
\[
\begin{aligned}
|\partial_xN_\varphi(t,x)-\partial_xN_\varphi(t,x')|&\leq|x-x'|\,\|\varphi\|_\infty\int_0^t\!\!\int_D\sup_{\xi\in D}|\partial_x^2\Delta_yG(t-s,\xi,y)|\,dy\,ds\\
&\leq\pi C_{2,2}\big(T+6T^{1/6}\big)\|\varphi\|_\infty|x-x'| .
\end{aligned}
\]
\end{proof}

\begin{theorem}[H\"older regularity]\label{thm:regularity}
Let $u_0\in H^1(D)$ and let $u$ be the solution of Theorem~\ref{thm:global}. Almost surely:
\begin{enumerate}
\item $u$ is $\frac1{12}$-H\"older continuous in $t$, uniformly in $x$, and $\frac12$-H\"older continuous in $x$, uniformly in $t$, on $[0,T]\times D$;
\item for every $t_0\in(0,T)$ and $\gamma_1<\frac14$, $u$ is $\gamma_1$-H\"older in $t$ and Lipschitz in $x$ on $[t_0,T]\times D$;
\item for every $t\in(0,T]$, $u(t,\cdot)\in C^1(D)$, $\partial_xu(t,0)=\partial_xu(t,\pi)=0$, and for every $t_0\in(0,T)$ and $\gamma<\frac12$, $\partial_xu$ is $\gamma$-H\"older in $x$, uniformly in $t\in[t_0,T]$.
\end{enumerate}
\end{theorem}

\begin{proof}
We use \eqref{eq:mild} in the form $u=S(\cdot)u_0+N_\varphi+Z$ with $\varphi=f(u)$. By \eqref{eq:globalbound}, $\|\varphi\|_\infty\leq B_M<\infty$ with $M=C'_0+K_Z$.

(1) By Lemma~\ref{lem:conv}, $N_\varphi$ is $\frac12$-H\"older in $t$ and Lipschitz in $x$, with constants $CB_M$. By Lemma~\ref{lem:Zholder}, $Z$ is $\gamma_1$-H\"older in $t$ for every $\gamma_1<\frac14$, in particular for $\gamma_1=\frac1{12}$, and by Lemma~\ref{lem:Zx} it is Lipschitz in $x$ with constant $K_Z$. With Lemma~\ref{lem:initialreg}, each of the three terms is $\frac1{12}$-H\"older in $t$ and $\frac12$-H\"older in $x$, uniformly in the other variable.

(2) On $[t_0,T]\times D$, $S(\cdot)u_0$ is Lipschitz in both variables by Lemma~\ref{lem:initialreg}, $N_\varphi$ is $\frac12$-H\"older in $t$ and Lipschitz in $x$, and $Z$ is $\gamma_1$-H\"older in $t$ and Lipschitz in $x$. The exponents in time are thus limited by $Z$.

(3) The $C^1$ property and the boundary values follow from Lemma~\ref{lem:initialreg}, Lemma~\ref{lem:convC1} and Lemma~\ref{lem:Zx}. For the H\"older continuity of $\partial_xu=\partial_xS(\cdot)u_0+\partial_xN_\varphi+\partial_xZ$ on $[t_0,T]\times D$: the first term is smooth by Lemma~\ref{lem:initialreg}, the second is Lipschitz in $x$ by Lemma~\ref{lem:convC1}, and $\partial_xZ=Y$ is $\gamma$-H\"older in $x$ for every $\gamma<\frac12$ by Step~1 of the proof of Lemma~\ref{lem:Zx} and Lemma~\ref{lem:kolmogorov}.
\end{proof}

\begin{remark}
The exponent $\frac1{12}$ near $t=0$ comes only from the initial condition. If $\|u_0\|^2_{H^s}=\sum_k(1+k^2)^s|\hat u_{0,k}|^2<\infty$ for some $s\in[1,2]$, the computation of Lemma~\ref{lem:initialreg} gives the exponent $\frac{2s-1}{12}$ for the initial term, and therefore every exponent $\gamma_1<\frac14$ on $[0,T]$ as soon as $s=2$. The exponent $\frac14$ comes from the spectral sum $\sum_kk^2\min(h,\lambda_k^{-1})\approx h^{1/2}$ in the proof of Lemma~\ref{lem:Zholder}: the factor $k^2$ produced by the conservative noise is only partly compensated by $\lambda_k\approx k^6$. For the stochastic Cahn-Hilliard equation, the corresponding sums are $\sum_kk^2\min(h,k^{-4})\approx h^{1/4}$ with conservative noise and $\sum_k\min(h,k^{-4})\approx h^{3/4}$ with additive white noise, which give the time exponents $\frac18$ and $\frac38$.
\end{remark}
\section{Malliavin differentiability and existence of a density}\label{sec:density}

\subsection{Malliavin calculus}

Let $\cH=L^2([0,T]\times D)$. The map $h\mapsto W(h)=\int_0^T\!\int_Dh(s,y)\,W(ds,dy)$ is an isonormal Gaussian process on $\cH$, and we use the associated Malliavin derivative $D$, the space $\D^{1,2}$ with the norm $\|F\|^2_{1,2}=\E|F|^2+\E\|DF\|^2_{\cH}$, and the following facts from Nualart \cite{Nu06}:
\begin{enumerate}
\item[(M1)] (closedness) if $F_n\in\D^{1,2}$, $F_n\to F$ in $L^2(\Omega)$ and $DF_n\to X$ in $L^2(\Omega;\cH)$, then $F\in\D^{1,2}$ and $DF=X$ (\cite[Section~1.2]{Nu06});
\item[(M2)] (chain rule) if $F\in\D^{1,2}$ and $\phi\in C^1(\R)$ has a bounded derivative, then $\phi(F)\in\D^{1,2}$ and $D\phi(F)=\phi'(F)DF$ (\cite[Proposition~1.2.3]{Nu06});
\item[(M3)] (localization) $F\in\D^{1,2}_{\mathrm{loc}}$ if there are events $\Omega_n\uparrow\Omega$ almost surely and $F_n\in\D^{1,2}$ with $F=F_n$ on $\Omega_n$; then $DF:=DF_n$ on $\Omega_n$ is well defined (\cite[Section~1.3.5]{Nu06});
\item[(M4)] (Bouleau-Hirsch criterion) if $F\in\D^{1,1}_{\mathrm{loc}}$ and $\|DF\|_{\cH}>0$ almost surely, then the law of $F$ is absolutely continuous with respect to the Lebesgue measure (\cite[Theorem~2.1.3]{Nu06}). Note that $\D^{1,2}_{\mathrm{loc}}\subset\D^{1,1}_{\mathrm{loc}}$.
\end{enumerate}
For $h\in\cH$ deterministic, $W(h)\in\D^{1,2}$ and $DW(h)=h$. In particular, by \eqref{eq:Z}, $Z(t,x)\in\D^{1,2}$ and
\begin{equation}\label{eq:DZ}
D_{r,z}Z(t,x)=h_{t,x}(r,z):=-\sigma\,\partial_zG(t-r,x,z)\,\mathbf 1_{\{r<t\}},\qquad \|h_{t,x}\|^2_{\cH}=\E|Z(t,x)|^2\leq C ,
\end{equation}

\begin{lemma}[Derivative of an integral]\label{lem:Dint}
Let $K\in L^1([0,T]\times D)$ be deterministic, and let $Y=\{Y(s,y)\}$ be a jointly measurable random field such that $Y(s,y)\in\D^{1,2}$ for every $(s,y)$, $(s,y,\omega)\mapsto DY(s,y)$ is jointly measurable, and $\sup_{s,y}\|Y(s,y)\|_{1,2}<\infty$. Then $F=\int_0^T\!\int_DK(s,y)Y(s,y)\,dy\,ds\in\D^{1,2}$ and $DF=\int_0^T\!\int_DK(s,y)\,DY(s,y)\,dy\,ds$, the integral being a Bochner integral in $L^2(\Omega;\cH)$.
\end{lemma}

\begin{proof}
$\D^{1,2}$ is a separable Hilbert space for $\|\cdot\|_{1,2}$, and the assumptions make $(s,y)\mapsto Y(s,y)\in\D^{1,2}$ strongly measurable (Pettis theorem) and Bochner integrable with respect to $|K(s,y)|\,dy\,ds$. The maps $\D^{1,2}\to L^2(\Omega)$, $F\mapsto F$, and $D:\D^{1,2}\to L^2(\Omega;\cH)$ are bounded and linear, hence commute with the Bochner integral (Hille's theorem). Finally, by Fubini's theorem, the Bochner integral in $L^2(\Omega)$ coincides almost surely with the pathwise integral $F$.
\end{proof}

\subsection{Differentiability of the truncated solution}

Recall the truncation $f_R=f\circ\rho_R$ of Section~\ref{sec:local}: $f_R\in C^\infty$ with $|f_R|\leq B_R$, $|f_R'|\leq L_R$, and $f_R'$ is Lipschitz with some constant $L'_R$.

\begin{proposition}\label{prop:D12}
Let $R>\|u_0\|_\infty$. For every $(t,x)\in[0,T]\times D$, $u_R(t,x)\in\D^{1,2}$. The random field $X=Du_R$ is the unique jointly measurable $\cH$-valued random field such that $\sup_{t,x}\|X(t,x)\|_{\cH}<\infty$ almost surely and
\begin{equation}\label{eq:Du}
X(t,x)=h_{t,x}+\int_0^t\!\!\int_D\Delta_yG(t-s,x,y)\,f_R'(u_R(s,y))\,X(s,y)\,dy\,ds,
\end{equation}
the integral being a Bochner integral in $\cH$, computed for each $\omega$. In particular $D_{r,z}u_R(t,x)=0$ for $r>t$, and there is a deterministic constant $C_R$ such that $\sup_{t,x}\|Du_R(t,x)\|_{\cH}\leq C_R$ almost surely.
\end{proposition}

\begin{proof}
\emph{Step 1: the linear equation.} For each $\omega$, let $E_\cH$ be the space of measurable maps $X:[0,T]\times D\to\cH$ with $\|X\|_{\beta}:=\sup_{t,x}e^{-\beta t}\|X(t,x)\|_{\cH}<\infty$, and let $\Psi(X)$ be the right-hand side of \eqref{eq:Du}. By Minkowski's inequality for Bochner integrals and \eqref{eq:conv1},
\[
e^{-\beta t}\|\Psi(X)(t,x)-\Psi(X')(t,x)\|_{\cH}\leq CL_R\int_0^t\big(1+(t-s)^{-1/2}\big)e^{-\beta(t-s)}ds\;\|X-X'\|_\beta,
\]
and $\|h_{t,x}\|_\cH\leq C$ by \eqref{eq:DZ}. As in Proposition~\ref{prop:trunc}, for $\beta=\beta_R$ deterministic, $\Psi$ is a contraction of $E_\cH$ with constant $\frac12$, and \eqref{eq:Du} has a unique solution $X$, with $\|X\|_{\beta_R}\leq2\sup\|h_{t,x}\|_\cH$. Hence $\sup_{t,x}\|X(t,x)\|_\cH\leq C_R:=2Ce^{\beta_RT}$. The Picard iterates $X^0=h$, $X^{m+1}=\Psi(X^m)$ vanish for $r>t$, and so does $X$. Each $X^m$ is jointly measurable in $(t,x,\omega)$ as an $\cH$-valued function: $h$ is deterministic and continuous in $(t,x)$ with values in $\cH$, and if $X^m$ is jointly measurable and bounded, then $\Psi(X^m)(t,x)$ is a Bochner integral of the jointly measurable integrand $(s,y,\omega)\mapsto\Delta_yG(t-s,x,y)f_R'(u_R(s,y))X^m(s,y)$, hence jointly measurable by Fubini's theorem for Bochner integrals. Since $X^m\to X$ uniformly in $(t,x)$ for every $\omega$, $X$ is jointly measurable.

\emph{Step 2: differentiability of the Picard iterates.} Let $u^0=0$ and $u^{n+1}=S(\cdot)u_0+N_{f_R(u^n)}+Z$ as in Proposition~\ref{prop:trunc}. We show by induction that $u^n(t,x)\in\D^{1,2}$ for all $(t,x)$, with
\begin{equation}\label{eq:Dun}
Du^{n+1}(t,x)=h_{t,x}+\int_0^t\!\!\int_D\Delta_yG(t-s,x,y)\,f_R'(u^n(s,y))\,Du^n(s,y)\,dy\,ds,
\end{equation}
and $\sup_{t,x}\|Du^n(t,x)\|_\cH\leq C_R$ almost surely. This is clear for $n=0$. If it holds for $n$, then by (M2) $f_R(u^n(s,y))\in\D^{1,2}$ with $Df_R(u^n)=f_R'(u^n)Du^n$, $|f_R(u^n)|\leq B_R$ and $\|Df_R(u^n)\|_\cH\leq L_RC_R$. Lemma~\ref{lem:Dint} with $K=\Delta_yG(t-\cdot,x,\cdot)\mathbf 1_{[0,t]}$, which is integrable by Lemma~\ref{lem:green}, gives $N_{f_R(u^n)}(t,x)\in\D^{1,2}$ and its derivative; since $S(t)u_0(x)$ is deterministic and $Z(t,x)\in\D^{1,2}$, \eqref{eq:Dun} follows. The bound on $Du^{n+1}$ follows from \eqref{eq:Dun} in the same way as in Step~1: with $\psi_n=\sup_{t,x}e^{-\beta_Rt}\|Du^n(t,x)\|_\cH$ and $C=\sup_{t,x}\|h_{t,x}\|_\cH$, we get $\psi_{n+1}\leq C+\frac12\psi_n$, so that $\psi_n\leq2C$ implies $\psi_{n+1}\leq2C$, and hence $\|Du^{n+1}(t,x)\|_\cH\leq2Ce^{\beta_RT}=C_R$.

\emph{Step 3: convergence.} Let $\epsilon_n=\sup_{t,x}|u^n(t,x)-u_R(t,x)|$, which satisfies $\epsilon_n\to0$ almost surely and $\E\epsilon_n^2\to0$ by \eqref{eq:picard}, since $\E\sup|Z|^2<\infty$ (Lemma~\ref{lem:Zx}). Subtracting \eqref{eq:Du} from \eqref{eq:Dun},
\begin{multline*}
Du^{n+1}(t,x)-X(t,x)\\
=\int_0^t\!\!\int_D\Delta_yG(t-s,x,y)\Big[f_R'(u^n)\big(Du^n-X\big)+\big(f_R'(u^n)-f_R'(u_R)\big)X\Big](s,y)\,dy\,ds,
\end{multline*}
and $|f_R'(u^n)-f_R'(u_R)|\leq L'_R\epsilon_n$. Let $E_n=\sup_{t,x}e^{-\beta_Rt}\|Du^n(t,x)-X(t,x)\|_\cH$, which is at most $2C_R$ by Step~2 and Step~1. By Minkowski's inequality in $\cH$, $|f_R'|\leq L_R$, $\|X\|_\cH\leq C_R$ and \eqref{eq:conv1}, for every $(t,x)$,
\begin{align*}
e^{-\beta_Rt}\|Du^{n+1}(t,x)-X(t,x)\|_\cH
&\leq e^{-\beta_Rt}\int_0^t\!\!\int_D|\Delta_yG(t-s,x,y)|\\
&\qquad\qquad\times\Big(L_R\|Du^n(s,y)-X(s,y)\|_\cH+L'_R\epsilon_nC_R\Big)dy\,ds\\
&\leq L_R\,E_n\int_0^t\!\!\int_D|\Delta_yG(t-s,x,y)|\,e^{-\beta_R(t-s)}\,dy\,ds\\
&\qquad+L'_RC_R\,\epsilon_n\int_0^t\!\!\int_D|\Delta_yG(t-s,x,y)|\,dy\,ds\\
&\leq\pi C_{0,2}L_R\big(\beta_R^{-1}+\Gamma(\tfrac12)\beta_R^{-1/2}\big)E_n\\
&\qquad+\pi C_{0,2}\big(T+2T^{1/2}\big)L'_RC_R\,\epsilon_n\\
&\leq\tfrac12E_n+C_5\,\epsilon_n,\qquad C_5:=\pi C_{0,2}\big(T+2T^{1/2}\big)L'_RC_R ,
\end{align*}
by the choice of $\beta_R$ in Step~1. Taking the supremum over $(t,x)$ and iterating,
\[
E_n\leq2^{-n}E_0+C_5\sum_{m=0}^{n-1}2^{-(n-1-m)}\epsilon_m\leq2^{1-n}C_R+C_5\sum_{m=0}^{n-1}2^{-(n-1-m)}\epsilon_m ,
\]
and the last sum tends to $0$ since $\epsilon_m\to0$ (a discrete convolution of a summable sequence with a sequence tending to zero). Hence $E_n\to0$ almost surely. By dominated convergence, $\E\|Du^n(t,x)-X(t,x)\|^2_\cH\to0$, and $u^n(t,x)\to u_R(t,x)$ in $L^2(\Omega)$. By (M1), $u_R(t,x)\in\D^{1,2}$ and $Du_R(t,x)=X(t,x)$.
\end{proof}

\begin{corollary}\label{cor:loc}
For $R>\|u_0\|_\infty$, let $\Omega_R=\big\{\sup_{t,x}|u(t,x)|<R\big\}$. Then $u=u_R$ on $[0,T]\times D$ on the event $\Omega_R$, the events $\Omega_R$ increase to an event of probability one, and for every $(t,x)$, $u(t,x)\in\D^{1,2}_{\mathrm{loc}}$ with $Du(t,x)=Du_R(t,x)$ on $\Omega_R$.
\end{corollary}

\begin{proof}
On $\Omega_R$, $\tau_R=T$: otherwise, by Theorem~\ref{thm:local}, $\tau_R<T$ and $\|u(\tau_R)\|_\infty=\|u_R(\tau_R)\|_\infty=R$. Hence $u=u_R$ on $[0,T]\times D$ on $\Omega_R$. By \eqref{eq:globalbound}, $\bigcup_R\Omega_R$ has probability one. The last claim follows from (M3) and Proposition~\ref{prop:D12}.
\end{proof}

\subsection{Non-degeneracy}

\begin{lemma}\label{lem:lower}
For every $x\in D$ there exist $c(x)>0$ and $\delta_0(x)\in(0,1]$ such that, for all $\delta\in(0,\delta_0(x)]$,
\[
c(x)\,\delta^{1/2}\leq\int_0^\delta\!\!\int_D|\partial_zG(r,x,z)|^2\,dz\,dr\leq C\,\delta^{1/2},
\]
with $C$ independent of $x$.
\end{lemma}

\begin{proof}
By \eqref{eq:normG}, the integral equals $e_1(x)^2\delta+\sum_{k\geq2}k^2e_k(x)^2\frac{1-e^{-2\lambda_k\delta}}{2\lambda_k}$.

\emph{Upper bound.} Since $\frac{1-e^{-2\lambda\delta}}{2\lambda}\leq\min(\delta,\frac1{2\lambda})$ and $\lambda_k\geq\frac9{16}k^6$, Lemma~\ref{lem:sums}(b) with $m=2$ bounds the sum by $C\delta^{1/2}$, and $e_1(x)^2\delta\leq\delta\leq\delta^{1/2}$.

\emph{Lower bound.} Let $N=\lfloor\delta^{-1/6}\rfloor$, so that $N\geq\frac12\delta^{-1/6}$ for $\delta\leq1$. For $2\leq k\leq N$, $\lambda_k\leq k^6\leq N^6\leq\delta^{-1}$, so $a:=\lambda_k\delta\leq1$, and since $a\mapsto(1-e^{-2a})/(2a)$ is decreasing, $\frac{1-e^{-2\lambda_k\delta}}{2\lambda_k}\geq c_*\delta$ with $c_*=\frac{1-e^{-2}}2$. Hence the integral is at least $\frac2\pi c_*\delta\sum_{k=2}^Nk^2\cos^2(kx)$. We claim that
\begin{equation}\label{eq:abel}
\sum_{k=2}^Nk^2\cos^2(kx)\geq\frac{N^3}{12}\qquad\text{for }N\geq N_0(x),
\end{equation}
where $N_0(x)=3$ if $x\in\{0,\pi\}$ and $N_0(x)=\lceil24/\sin x\rceil$ otherwise. If $x\in\{0,\pi\}$, $\cos^2(kx)=1$ and $\sum_{k=2}^Nk^2\geq\frac{N^3}3-1\geq\frac{N^3}{12}$ for $N\geq3$. If $x\in(0,\pi)$, write $\cos^2(kx)=\frac12(1+\cos(2kx))$. The partial sums $S_m=\sum_{k=1}^m\cos(2kx)=\frac{\sin(mx)\cos((m+1)x)}{\sin x}$ satisfy $|S_m|\leq1/\sin x$, and Abel's summation gives
\[
\Big|\sum_{k=1}^Nk^2\cos(2kx)\Big|=\Big|N^2S_N-\sum_{k=1}^{N-1}\big((k+1)^2-k^2\big)S_k\Big|\leq\frac{N^2+(N^2-1)}{\sin x}\leq\frac{2N^2}{\sin x}.
\]
Hence $\sum_{k=2}^Nk^2\cos^2(kx)\geq\frac12\big(\frac{N^3}3-1\big)-\frac12\big(\frac{2N^2}{\sin x}+1\big)=\frac{N^3}6-\frac{N^2}{\sin x}-1$, which is at least $\frac{N^3}{12}$ as soon as $\frac{N^3}{12}\geq\frac{N^2}{\sin x}+1$, in particular for $N\geq24/\sin x$. This proves \eqref{eq:abel}. For $\delta\leq\delta_0(x):=(2N_0(x))^{-6}$ we have $N\geq N_0(x)$, and the integral is at least $\frac2\pi c_*\delta\frac{N^3}{12}\geq\frac{c_*}{48\pi}\delta^{1/2}$.
\end{proof}

\begin{proposition}\label{prop:nondeg}
Assume $\sigma\neq0$ and $R>\|u_0\|_\infty$. For every $t\in(0,T]$ and $x\in D$, $\|Du_R(t,x)\|_{\cH}>0$ almost surely.
\end{proposition}

\begin{proof}
Fix $t\in(0,T]$, $x\in D$ and $\delta\in(0,t\wedge\delta_0(x)]$. For $(s,y)\in[t-\delta,t]\times D$ let
\[
Q(s,y)=\Big(\int_{t-\delta}^s\!\!\int_D|D_{r,z}u_R(s,y)|^2\,dz\,dr\Big)^{1/2}\leq\|Du_R(s,y)\|_\cH\leq C_R .
\]

\emph{Step 1: a bound on $Q$.} Let $(s,y)\in[t-\delta,t]\times D$ and $r\in[t-\delta,s]$. Since $D_{r,z}u_R(s',y')=0$ for $s'<r$, \eqref{eq:Du} at $(s,y)$ reads
\begin{multline*}
D_{r,z}u_R(s,y)=-\sigma\,\partial_zG(s-r,y,z)\\
+\int_{t-\delta}^s\!\!\int_D\Delta_{y'}G(s-s',y,y')\,f_R'(u_R(s',y'))\,\mathbf 1_{\{r<s'\}}D_{r,z}u_R(s',y')\,dy'\,ds'.
\end{multline*}
We take the norm of $L^2([t-\delta,s]\times D,dr\,dz)$. For the first term, the change of variables $r'=s-r\in[0,s-t+\delta]\subset[0,\delta]$ and Lemma~\ref{lem:lower} give
\begin{align*}
|\sigma|\Big(\int_{t-\delta}^s\!\!\int_D|\partial_zG(s-r,y,z)|^2dz\,dr\Big)^{1/2}
&\leq|\sigma|\Big(\int_0^\delta\!\!\int_D|\partial_zG(r',y,z)|^2dz\,dr'\Big)^{1/2}\\
&\leq|\sigma|C^{1/2}\delta^{1/4}.
\end{align*}
For the second term, we use Minkowski's integral inequality, $|f_R'|\leq L_R$, and the fact that the $L^2([t-\delta,s]\times D)$ norm of $(r,z)\mapsto\mathbf 1_{\{r<s'\}}D_{r,z}u_R(s',y')$ equals $Q(s',y')$, since $D_{r,z}u_R(s',y')=0$ for $r>s'$. This gives the bound
\[
L_R\int_{t-\delta}^s\!\!\int_D|\Delta_{y'}G(s-s',y,y')|\,Q(s',y')\,dy'\,ds' .
\]
Hence, with $k(r,y,y')=|\Delta_{y'}G(r,y,y')|/(\pi C_{0,2})$, which satisfies $\int_Dk(r,y,y')dy'\leq1+r^{-1/2}$ by Lemma~\ref{lem:green},
\[
Q(s,y)\leq|\sigma|C^{1/2}\delta^{1/4}+\pi C_{0,2}L_R\int_{t-\delta}^s\!\!\int_Dk(s-s',y,y')\,Q(s',y')\,dy'\,ds' .
\]
Lemma~\ref{lem:gronwall} on $I=[t-\delta,t]$, with $a=|\sigma|C^{1/2}\delta^{1/4}$ and $b=\pi C_{0,2}L_R$, gives
\[
\sup_{(s,y)\in[t-\delta,t]\times D}Q(s,y)\leq2e^{\gamma(b)\delta}|\sigma|C^{1/2}\delta^{1/4}\leq C'_R\,\delta^{1/4},\qquad C'_R:=2e^{\gamma(b)T}|\sigma|C^{1/2}.
\]

\emph{Step 2: the remainder.} For $r\in[t-\delta,t]$ write $D_{r,z}u_R(t,x)=-\sigma\,\partial_zG(t-r,x,z)+I_{r,z}$, with
\[
I_{r,z}=\int_{t-\delta}^t\!\!\int_D\Delta_yG(t-s,x,y)\,f_R'(u_R(s,y))\,\mathbf 1_{\{r<s\}}D_{r,z}u_R(s,y)\,dy\,ds .
\]
As in Step~1, by Minkowski's inequality, Lemma~\ref{lem:green} and Step~1,
\begin{align*}
\Big(\int_{t-\delta}^t\!\!\int_D|I_{r,z}|^2\,dz\,dr\Big)^{1/2}
&\leq L_R\int_{t-\delta}^t\!\!\int_D|\Delta_yG(t-s,x,y)|\,Q(s,y)\,dy\,ds\\
&\leq\pi C_{0,2}L_R\,C'_R\,\delta^{1/4}\int_{t-\delta}^t\big(1+(t-s)^{-1/2}\big)ds\\
&=\pi C_{0,2}L_RC'_R\,\delta^{1/4}\big(\delta+2\delta^{1/2}\big)
\leq3\pi C_{0,2}L_RC'_R\,\delta^{3/4}=:C''_R\,\delta^{3/4},
\end{align*}
where we used $\delta\leq\delta^{1/2}$ for $\delta\leq1$.

\emph{Step 3: conclusion.} Using the elementary inequality $|a+b|^2\geq\frac12|a|^2-|b|^2$ and Lemma~\ref{lem:lower},
\begin{align*}
\|Du_R(t,x)\|^2_\cH
&\geq\int_{t-\delta}^t\!\!\int_D|D_{r,z}u_R(t,x)|^2\,dz\,dr
=\int_{t-\delta}^t\!\!\int_D\big|-\sigma\partial_zG(t-r,x,z)+I_{r,z}\big|^2\,dz\,dr\\
&\geq\frac{\sigma^2}2\int_0^\delta\!\!\int_D|\partial_zG(r,x,z)|^2\,dz\,dr-\int_{t-\delta}^t\!\!\int_D|I_{r,z}|^2\,dz\,dr\\
&\geq\frac{\sigma^2c(x)}2\,\delta^{1/2}-(C''_R)^2\,\delta^{3/2}
=\delta^{1/2}\Big(\frac{\sigma^2c(x)}2-(C''_R)^2\,\delta\Big).
\end{align*}
Choosing $\delta\leq\min\big(t,\delta_0(x),\sigma^2c(x)/(4(C''_R)^2)\big)$, the right-hand side is at least $\frac{\sigma^2c(x)}4\delta^{1/2}>0$. Since $c(x)$, $C''_R$ and $\delta$ are deterministic and the inequalities hold almost surely, $\|Du_R(t,x)\|_\cH>0$ almost surely.
\end{proof}

\subsection{Main result}

\begin{theorem}[Existence of a density]\label{thm:density}
Let $u_0\in H^1(D)$ and $\sigma\neq0$. For every $t\in(0,T]$ and $x\in D$, the law of $u(t,x)$ is absolutely continuous with respect to the Lebesgue measure on $\R$.
\end{theorem}

\begin{proof}
By Corollary~\ref{cor:loc}, $u(t,x)\in\D^{1,2}_{\mathrm{loc}}$, localized by $(\Omega_R,u_R(t,x))$, and $\|Du(t,x)\|_\cH=\|Du_R(t,x)\|_\cH$ on $\Omega_R$. By Proposition~\ref{prop:nondeg}, $\|Du(t,x)\|_\cH>0$ almost surely on each $\Omega_R$, hence almost surely. The criterion (M4) concludes.
\end{proof}

\begin{remark}
\begin{enumerate}
\item Theorem~\ref{thm:density} holds up to the boundary, $x\in\{0,\pi\}$ included, because the lower bound of Lemma~\ref{lem:lower} also covers the boundary points, with $\delta_0(0)=\delta_0(\pi)=6^{-6}$. For the stochastic Cahn-Hilliard equation, the existence of a density is stated in \cite{CW01} for $x$ in the interior of the domain.
\item For $t=0$ the law of $u(0,x)=u_0(x)$ is a Dirac mass, and the conservation of the mean shows that $\int_Du(t,y)\,dy$ has no density either. Theorem~\ref{thm:density} concerns the pointwise values $u(t,x)$, for which there is no such constraint.
\item The smoothness of the density would require $u(t,x)\in\D^\infty$ and $\E\|Du(t,x)\|_\cH^{-p}<\infty$ for all $p$. The moment bound \eqref{eq:moments} is a first step in this direction, but the non-degeneracy estimate of Proposition~\ref{prop:nondeg} is qualitative, its constants depending on the localization level $R$. This question is left open.
\end{enumerate}
\end{remark}
\section{Numerical illustrations}\label{sec:numerics}

In this section we illustrate the results of Sections~\ref{sec:global}, \ref{sec:regularity} and~\ref{sec:density} on numerical simulations. The scheme below is not analysed in this paper; its convergence is a natural question for future work, for which the results above provide the necessary regularity. All Monte Carlo estimates are given with their $95\%$ confidence intervals, or with the corresponding relative accuracy.

\subsection{The scheme}

We use the spectral Galerkin approximation of Section~\ref{sec:global} with $N$ modes, combined with an exponential Euler step in time. Writing $\xi_k(t)=\langle u(t),e_k\rangle$, the projected equation reads
\[
d\xi_k=\big(-\lambda_k\xi_k-k^2\langle f(u),e_k\rangle\big)\,dt+\sigma k\,dB_k,\qquad 1\leq k\leq N,\qquad d\xi_0=0,
\]
where the noise term follows from \eqref{eq:Zseries}. For a time step $\Delta>0$ we freeze the nonlinearity on each step and solve the linear part and the noise exactly:
\begin{equation}\label{eq:scheme}
\xi_k^{n+1}=e^{-\lambda_k\Delta}\xi_k^n-k^2\,\frac{1-e^{-\lambda_k\Delta}}{\lambda_k}\,\big\langle f(u^n),e_k\big\rangle+\sigma k\Big(\frac{1-e^{-2\lambda_k\Delta}}{2\lambda_k}\Big)^{1/2}\eta_k^n,
\end{equation}
where $(\eta_k^n)$ are independent standard Gaussian variables and the fractions are replaced by $\Delta$ when $\lambda_k=0$, that is for $k=1$. The coefficients $\langle f(u^n),e_k\rangle$ are computed with a discrete cosine transform on $2N$ midpoints. The midpoint rule on $2N$ points integrates exactly the trigonometric polynomials $\cos(jx)$ with $j<4N$, so that $\langle(u^n)^3,e_k\rangle$ is computed exactly for $k<N$; only the last coefficient $k=N$ is affected by aliasing. We note three features of \eqref{eq:scheme}.
\begin{enumerate}
\item The mean $\xi_0$ is exactly constant, as in Proposition~\ref{prop:mass}.
\item For the retained modes $k\leq N$, the stochastic convolution of the linear equation is sampled exactly in law at the grid times, so that no time discretization error is made on it.
\item For $k\geq2$, the coefficient in front of the nonlinearity is at most $k^2/\lambda_k\leq\frac{16}9k^{-4}\leq\frac19$, whatever $\Delta$. Hence the stiffness of the high modes ($\lambda_k\approx k^6$) does not constrain the time step. Only the neutral mode $e_1$ is treated by an explicit Euler step, which requires $\Delta\,|f'(u)|$ to be small; this is the case for the values of $\Delta$ used below, since the recorded values of $|u|$ never exceed $1.16$ in our simulations, so that $\Delta\,|f'(u)|\leq\Delta(3\cdot1.16^2+\varepsilon)<0.04$ for $\Delta\leq4\cdot10^{-3}$.
\end{enumerate}
In all experiments $\varepsilon=0.5$, $\sigma=0.2$ and $u_0=0$. The Python code, with the random seeds, is available from the author.

\subsection{Choice of the resolution}

Table~\ref{tab:resolution} reports two quantities, $\E[\xi_1(6)^2]$ and $\E[u(3,0)^2]$, for several values of $N$ and $\Delta$, each estimated over $2000$ independent paths. No dependence on $N$ or on $\Delta$ can be detected at the Monte Carlo accuracy, which is between $2.5\%$ and $4.5\%$ here. This justifies the resolutions used below: $N=64$ and $\Delta=2\cdot10^{-3}$ for ensembles, $N=128$ and $\Delta=10^{-3}$ for the single path of Figure~\ref{fig:symmetry}(a), and $N=256$ for the spatial increments of Section~\ref{sec:numreg}.

\begin{table}[H]
\centering
\begin{tabular}{rrcc}
\hline
$N$ & $\Delta$ & $\E[\xi_1(6)^2]$ & $\E[u(3,0)^2]$\\
\hline
$32$ & $2\cdot10^{-3}$ & $0.753\pm0.018$ & $0.220\pm0.009$\\
$64$ & $4\cdot10^{-3}$ & $0.749\pm0.018$ & $0.222\pm0.010$\\
$64$ & $2\cdot10^{-3}$ & $0.734\pm0.018$ & $0.218\pm0.009$\\
$64$ & $1\cdot10^{-3}$ & $0.744\pm0.018$ & $0.214\pm0.009$\\
$128$ & $2\cdot10^{-3}$ & $0.759\pm0.018$ & $0.219\pm0.009$\\
\hline
\end{tabular}
\caption{Dependence on the number of modes $N$ and on the time step $\Delta$, with $95\%$ confidence intervals over $2000$ paths. The differences are within the Monte Carlo uncertainty.}
\label{tab:resolution}
\end{table}

\subsection{Symmetry breaking triggered by the noise}

Since $f(0)=0$, the constant function $u\equiv0$ is a stationary solution of the deterministic equation. It is linearly unstable: the linearization at $0$ acts on $e_k$ by multiplication by $-\lambda_k+\varepsilon k^2$, which is positive only for $k=1$, because $\lambda_1=0$. Neglecting the higher modes, which are strongly damped ($\lambda_3=576$), the projection on $e_1$ satisfies
\[
\xi_1'=-\big\langle f(\xi_1e_1),e_1\big\rangle=\varepsilon\xi_1-\xi_1^3\int_De_1^4\,dx=\varepsilon\xi_1-\frac3{2\pi}\,\xi_1^3,
\]
whose stable equilibria are $\xi_1=\pm(2\pi\varepsilon/3)^{1/2}=\pm\sqrt{\pi/3}$ for $\varepsilon=0.5$. They correspond to the two periodic profiles $u\approx\pm(4\varepsilon/3)^{1/2}\cos x\approx\pm0.82\cos x$, the one dimensional analogues of the crystalline phase, while $u\equiv0$ plays the role of the homogeneous (liquid) phase.

\begin{figure}[H]
\centering
\includegraphics[width=\textwidth]{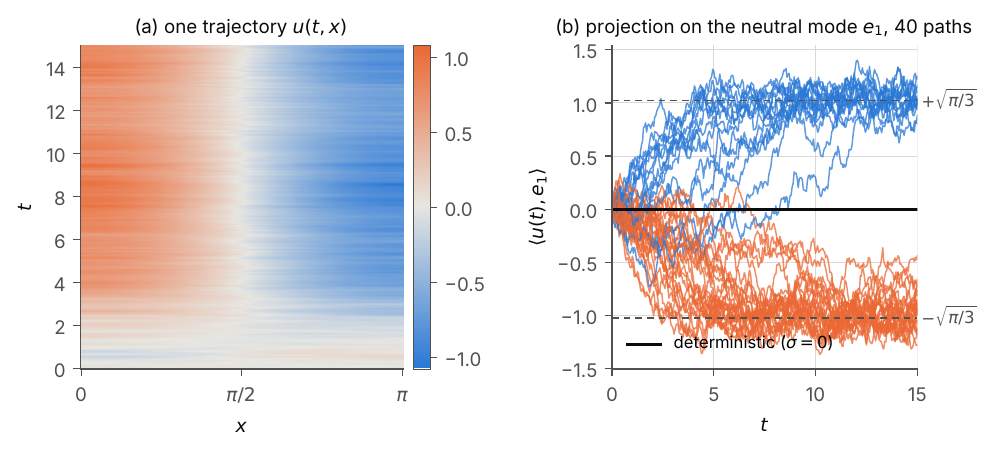}
\caption{Noise-induced symmetry breaking, $\varepsilon=0.5$, $\sigma=0.2$, $u_0=0$. (a) One path of $u(t,x)$, $N=128$, $\Delta=10^{-3}$. (b) Projection $\langle u(t),e_1\rangle$ for $40$ independent paths ($N=64$, $\Delta=2\cdot10^{-3}$), coloured according to the selected branch. The deterministic solution ($\sigma=0$) stays at $0$. Dashed lines: the stable equilibria $\pm\sqrt{\pi/3}$ of the one-mode reduction.}
\label{fig:symmetry}
\end{figure}

Figure~\ref{fig:symmetry} shows that the deterministic scheme started from $u_0=0$ stays exactly at $0$, whereas the noise drives every path away from the homogeneous state, towards one of the two profiles, chosen at random. Near $0$, the linearized dynamics of $\xi_1$ is an unstable Ornstein-Uhlenbeck process, $\xi_1(t)\approx\sigma\int_0^te^{\varepsilon(t-s)}dB_1(s)$, whose size is of order $\sigma e^{\varepsilon t}(2\varepsilon)^{-1/2}$; it reaches order one after a time of order $\varepsilon^{-1}\log\big((2\varepsilon)^{1/2}/\sigma\big)\approx3.2$. Over the $40$ paths of Figure~\ref{fig:symmetry}(b), the first time at which $|\xi_1|$ exceeds half of $\sqrt{\pi/3}$ has median $2.8$, and ranges from $0.85$ to $8.5$. This is the role that the noise plays in the physics literature, where it is needed to leave the homogeneous state. The mode which carries the transition is the neutral mode $e_1$ of Remark~\ref{rem:neutral}, which is not damped by the linear part. Over $4000$ independent paths, the empirical mean of $\xi_1(15)^2$ is $0.988$, slightly below the value $\pi/3\approx1.047$ of the one-mode reduction; the difference comes from the fluctuations and from the coupling with the higher modes. On $D=[0,\pi]$ the transition is supercritical, so that the experiment illustrates a noise-induced symmetry breaking rather than a metastable nucleation event; the latter requires larger domains, where several unstable modes coexist.

\subsection{H\"older regularity}\label{sec:numreg}

\begin{figure}[H]
\centering
\includegraphics[width=\textwidth]{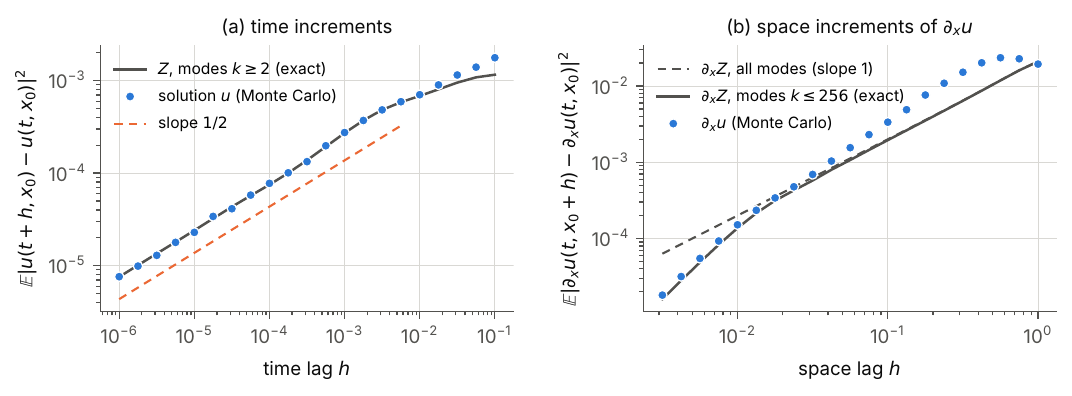}
\caption{Mean square increments at $x_0=\pi/3$, in logarithmic scales. (a) Time increments of $u$ over $2000$ paths ($N=64$), compared with the exact increments of the modes $k\geq2$ of the stochastic convolution and with the slope $\frac12$. (b) Space increments of $\partial_xu$ over $1000$ paths ($N=256$), compared with the exact increments of $\partial_xZ$, computed with all modes (slope $1$) and with the modes $k\leq256$.}
\label{fig:variograms}
\end{figure}

\paragraph{Time.} We start from $2000$ paths at time $t=15$, which are close to the stationary regime, and estimate $\E|u(t+h,x_0)-u(t,x_0)|^2$ at $x_0=\pi/3$ for $h\in[10^{-6},10^{-1}]$, integrating each increment with $16$ steps of size $h/16$. With $2000$ paths, the relative accuracy of each estimate is about $\pm6\%$ at the $95\%$ level. The results are shown in Figure~\ref{fig:variograms}(a), together with the exact mean square increments of the stationary part of the stochastic convolution, that is of its modes $k\geq2$, computed from \eqref{eq:Zseries}:
\[
\E|Z_{\geq2}(t+h,x_0)-Z_{\geq2}(t,x_0)|^2=\sigma^2\sum_{k\geq2}k^2e_k(x_0)^2\,\frac{1-e^{-\lambda_kh}}{\lambda_k}.
\]
A least squares fit on $h\in[10^{-6},10^{-3}]$ gives the slope $0.516$, in agreement with the exponent $\frac12$ of Lemma~\ref{lem:Zholder}, which corresponds to the H\"older exponents $\gamma_1<\frac14$. For $h\leq10^{-2}$, the empirical increments of $u$ agree with the exact increments of $Z_{\geq2}$ within the Monte Carlo uncertainty (the largest relative deviation is $5.4\%$). This confirms the mechanism of the proof of Theorem~\ref{thm:regularity}: at small time scales the regularity of $u$ is that of $Z$, the nonlinear term being smoother ($\frac12$-H\"older by Lemma~\ref{lem:conv}). For $h\geq2\cdot10^{-2}$ the empirical increments are larger than those of $Z_{\geq2}$, because of the neutral mode $e_1$, which is not included in the reference curve and whose dynamics is driven by the nonlinearity.

\paragraph{Space.} For the spatial regularity, the increments of $u$ itself are not informative: they are of order $h^2$ in mean square, but so are those of any smooth function, in particular of the truncated series. We therefore test the sharper statement of Theorem~\ref{thm:regularity}(3), namely the H\"older regularity of $\partial_xu$, which is governed by $\partial_xZ$, whose mean square increments are of order $h$ (Lemma~\ref{lem:Zx}). Since the variance of $\partial_xZ$ involves the slowly converging series $\sum_kk^{-2}$, this test requires a fine resolution: starting from $1000$ paths of the ensemble at $t=15$, we extend them to $N=256$ modes and let the new modes relax during a time $0.3$, which is much larger than their relaxation times $\lambda_k^{-1}\leq\lambda_{65}^{-1}<10^{-10}$. Figure~\ref{fig:variograms}(b) compares the empirical mean square increments of $\partial_xu$ at $x_0=\pi/3$ with two exact references: the increments of $\partial_xZ$ computed with all modes, which have slope $1$, and with the modes $k\leq256$. For $h\leq0.024$, the empirical increments agree with the truncated reference within $13\%$, the relative accuracy of the estimates being about $\pm9\%$ with $1000$ paths. Two effects prevent a clean slope $1$ from being observed: for $h\lesssim0.03$ the truncation at $N=256$ steepens the curve, as shown by the truncated reference, and for $h\gtrsim0.03$ the smooth part of the solution, mainly carried by the neutral mode, dominates. The agreement with the truncated reference shows, however, that at small spatial scales the regularity of $\partial_xu$ is that of $\partial_xZ$, as in the proof of Theorem~\ref{thm:regularity}(3).

\subsection{Densities}

\begin{figure}[H]
\centering
\includegraphics[width=\textwidth]{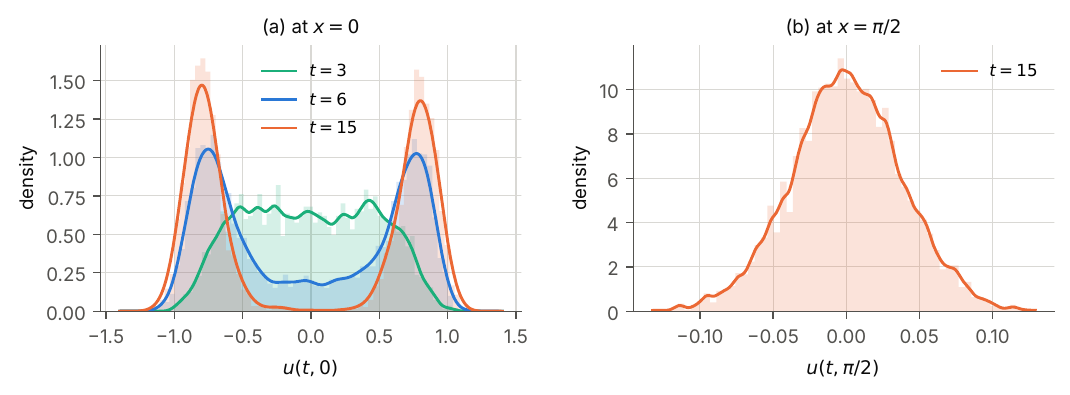}
\caption{Histograms and kernel density estimates of $u(t,x)$ over $4000$ paths. (a) At $x=0$ for $t=3,6,15$. (b) At $x=\pi/2$ for $t=15$.}
\label{fig:density}
\end{figure}

Figure~\ref{fig:density} shows histograms and kernel density estimates of $u(t,x)$ over $4000$ independent paths. At $x=0$, the law of $u(t,0)$ is broad and flat at $t=3$, when part of the paths has not yet left the homogeneous state, and becomes clearly bimodal at $t=6$ and $t=15$. At $t=15$ its two modes are located near $\pm0.78$, slightly inside the values $\pm0.82$ of the one-mode reduction, and $48.9\%$ of the paths have selected the positive branch, in agreement with the symmetry of the problem. At $x=\pi/2$ the profile $\cos x$ vanishes, the neutral mode does not contribute, and the law of $u(15,\pi/2)$ is unimodal and concentrated, with standard deviation $0.038$. In both cases the empirical distributions show no concentration at a point, consistently with Theorem~\ref{thm:density}. They also show that the density can be very different from a Gaussian one, which is why its existence requires the Malliavin calculus argument of Section~\ref{sec:density} rather than a Gaussian computation.

\section{Concluding remarks and open problems}\label{sec:conclusion}

We have shown that the one dimensional stochastic phase field crystal equation with conservative space-time white noise is globally well-posed, with moments of all orders, that its solution has the H\"older regularity dictated by the stochastic convolution, and that the law of $u(t,x)$ has a density for every $t>0$ and $x\in D$. The two features that distinguish the PFC equation from the Cahn-Hilliard equation, the second neutral mode and the order six of the operator, enter respectively through Lemmas~\ref{lem:interp} and~\ref{lem:Hminus1}, and through the spectral bounds of Section~\ref{sec:setting}. Several questions remain open.
\begin{enumerate}
\item \emph{Higher dimensions.} On $[0,\pi]^d$ the stochastic convolution has finite variance for $d\leq3$, but the a priori estimates of Section~\ref{sec:global} use the one dimensional embedding $H^{2/3}\subset L^\infty$ and the $C^1$ regularity of $Z$, which both fail in higher dimension.
\item \emph{Smoothness and positivity of the density.} The moment bounds of Theorem~\ref{thm:global} are a first step towards the smoothness of the density. The strict positivity of the density, proved by Cardon-Weber \cite{CW02} for the stochastic Cahn-Hilliard equation, is also open here.
\item \emph{Long time behaviour.} Since the mean is conserved, one expects an invariant measure for each value of the mean. The simulations of Section~\ref{sec:numerics} suggest that, for $\varepsilon>0$ and zero mean, it concentrates near the two profiles $\pm(4\varepsilon/3)^{1/2}\cos x$ when $\sigma$ is small.
\item \emph{Numerical analysis.} The convergence of the exponential scheme \eqref{eq:scheme}, with rates, is a natural continuation of this work.
\item \emph{Other noises.} The case of a noise which is fractional in time, and the non-Gaussian case of a Hermite noise, are left for future work.
\end{enumerate}

\section*{Data availability}
The Python code used for the simulations of Section~\ref{sec:numerics} is available from the author upon request.

\appendix
\section{Auxiliary results}\label{app:tools}

\begin{lemma}[Gronwall lemma with weakly singular kernel]\label{lem:gronwall}
Let $I=[t_0,t_1]$ with $t_1-t_0\leq T$, let $a,b\geq0$, and let $Q:I\times D\to[0,\infty)$ be bounded and measurable. Assume that for all $(s,y)\in I\times D$,
\[
Q(s,y)\leq a+b\int_{t_0}^s\!\!\int_D k(s-s',y,y')\,Q(s',y')\,dy'\,ds',
\]
where $k\geq0$ is measurable and $\int_Dk(r,y,y')\,dy'\leq1+r^{-1/2}$ for all $r>0$, $y\in D$. Let $\gamma=\gamma(b)>0$ be such that $b\big(\gamma^{-1}+\Gamma(\tfrac12)\gamma^{-1/2}\big)\leq\frac12$. Then $Q(s,y)\leq2a\,e^{\gamma(s-t_0)}$ on $I\times D$. The same holds for a bounded measurable $\varphi:I\to[0,\infty)$ with $\varphi(s)\leq a+b\int_{t_0}^s(1+(s-s')^{-1/2})\varphi(s')ds'$, by taking $Q(s,y)=\varphi(s)$ and $k=\pi^{-1}(1+r^{-1/2})$.
\end{lemma}

\begin{proof}
Let $M=\sup_{(s,y)\in I\times D}e^{-\gamma(s-t_0)}Q(s,y)$, which is finite since $Q$ is bounded. For $(s,y)\in I\times D$, using $Q(s',y')\leq Me^{\gamma(s'-t_0)}$ and the assumption on $k$,
\begin{align*}
e^{-\gamma(s-t_0)}Q(s,y)
&\leq e^{-\gamma(s-t_0)}a+b\int_{t_0}^s\!\!\int_Dk(s-s',y,y')\,e^{-\gamma(s-t_0)}Q(s',y')\,dy'\,ds'\\
&\leq a+bM\int_{t_0}^s\Big(\int_Dk(s-s',y,y')\,dy'\Big)e^{-\gamma(s-s')}\,ds'\\
&\leq a+bM\int_{t_0}^s\big(1+(s-s')^{-1/2}\big)e^{-\gamma(s-s')}\,ds'\\
&\leq a+bM\int_0^\infty\big(1+r^{-1/2}\big)e^{-\gamma r}\,dr
=a+bM\big(\gamma^{-1}+\Gamma(\tfrac12)\gamma^{-1/2}\big)\leq a+\frac M2 .
\end{align*}
Taking the supremum over $(s,y)$ gives $M\leq a+\frac M2$, hence $M\leq2a$, that is $Q(s,y)\leq2a\,e^{\gamma(s-t_0)}$. Only the measurability of $Q$ on $I\times D$ is used, not that of $s\mapsto\sup_yQ(s,y)$.
\end{proof}

\begin{lemma}[Spectral sums]\label{lem:sums}
\begin{enumerate}
\item[(a)] For $m\geq0$ and $c>0$ there is $C=C(m,c,T)$ such that $\sum_{k\geq2}k^me^{-ck^6t}\leq Ct^{-(m+1)/6}$ for $t\in(0,T]$.
\item[(b)] For $0\leq m<5$ there is $C=C(m,T)$ such that $\sum_{k\geq2}k^m\min(h,k^{-6})\leq Ch^{(5-m)/6}$ for $h\in(0,T]$.
\item[(c)] For $p>1$, $p\neq3$, there is $C=C(p)$ such that $\sum_{k\geq1}k^{-p}\min(k^2h^2,4)\leq Ch^{\min(p-1,2)}$ for $h\in(0,\pi]$.
\end{enumerate}
\end{lemma}

\begin{proof}
(a) The function $g(r)=r^me^{-cr^6t}$ is nondecreasing on $[0,r_*]$ and nonincreasing on $[r_*,\infty)$. For such a function, $\sum_{k\geq1}g(k)\leq\int_0^\infty g(r)\,dr+\sup_{r\geq0}g(r)$. The change of variables $r=t^{-1/6}\rho$ gives $\int_0^\infty g=t^{-(m+1)/6}\int_0^\infty\rho^me^{-c\rho^6}d\rho$ and $\sup g=t^{-m/6}\sup_\rho\rho^me^{-c\rho^6}$, and $t^{-m/6}\leq T^{1/6}t^{-(m+1)/6}$.

(b) Let $N=\lceil h^{-1/6}\rceil$, so that $h^{-1/6}\leq N\leq h^{-1/6}+1\leq(1+T^{1/6})h^{-1/6}$. Then $\sum_{2\leq k\leq N}k^mh\leq N^{m+1}h\leq Ch^{(5-m)/6}$, and $\sum_{k>N}k^{m-6}\leq\int_{N-1}^\infty r^{m-6}dr\leq CN^{m-5}\leq Ch^{(5-m)/6}$, where we used $m<5$ and $N\geq2$ when $h<1$ (for $h\geq1$ the bound is trivial).

(c) Let $N=\lfloor1/h\rfloor$. The terms $k\leq N$ contribute at most $h^2\sum_{k\leq N}k^{2-p}$, which is bounded by $Ch^2N^{3-p}\leq Ch^{p-1}$ if $p<3$, and by $Ch^2$ if $p>3$. The terms $k>N$ contribute at most $4\sum_{k>N}k^{-p}\leq C(N+1)^{1-p}\leq Ch^{p-1}$, and $h^{p-1}\leq\pi^{p-3}h^2$ if $p>3$.
\end{proof}

\begin{lemma}[Anisotropic Kolmogorov theorem for Gaussian fields]\label{lem:kolmogorov}
Let $X=\{X(t,x),(t,x)\in[0,T]\times D\}$ be a centred Gaussian random field such that
\begin{equation}\label{eq:kolhyp}
\E|X(t,x)-X(s,y)|^2\leq C\big(|t-s|^{2\alpha_1}+|x-y|^{2\alpha_2}\big),\qquad\alpha_1,\alpha_2\in(0,1].
\end{equation}
Then $X$ has a continuous version, and for all $\gamma_1<\alpha_1$ and $\gamma_2<\alpha_2$ there is a random variable $K$ with finite moments of all orders such that, for all $s,t\in[0,T]$ and $x,y\in D$,
\[
|X(t,x)-X(s,x)|\leq K|t-s|^{\gamma_1},\qquad|X(t,x)-X(t,y)|\leq K|x-y|^{\gamma_2}.
\]
\end{lemma}

\begin{proof}
\emph{Step 0: the Garsia-Rodemich-Rumsey inequality.} Let $I\subset\R$ be a bounded interval, $f:I\to\R$ continuous, $p\geq1$ and $\eta>0$. Theorem~1.1 of \cite{Wa86} in dimension one, with $\Psi(u)=|u|^p$ and $p(u)=u^{\eta+2/p}$, shows that if
\[
B=\int_I\!\!\int_I\frac{|f(r)-f(r')|^p}{|r-r'|^{\eta p+2}}\,dr\,dr'<\infty,
\]
then
\[|f(r)-f(r')|\leq c_{\eta,p}\,B^{1/p}\,|r-r'|^\eta\quad\text{for all }r,r'\in I .
\]

\emph{Step 1: continuity.} Since the increments of $X$ are Gaussian, \eqref{eq:kolhyp} gives, for every $p\geq1$ and $\alpha=\min(\alpha_1,\alpha_2)$,
\[
\E|X(t,x)-X(s,y)|^{2p}\leq C_p\big(|t-s|^{2\alpha_1}+|x-y|^{2\alpha_2}\big)^p\leq C'_p\,\big|(t,x)-(s,y)\big|^{2\alpha p}.
\]
Choosing $2\alpha p>2$, the Kolmogorov criterion \cite[Corollary~1.2]{Wa86} gives a continuous version of $X$, which we use from now on.

\emph{Step 2: spatial increments, uniformly in time.} Fix $x,y\in D$, $h=|x-y|$, and let $G(t)=X(t,x)-X(t,y)$. Writing $G(t)-G(s)$ in the two ways
\begin{align*}
G(t)-G(s)&=\big[X(t,x)-X(s,x)\big]-\big[X(t,y)-X(s,y)\big]\\
&=\big[X(t,x)-X(t,y)\big]-\big[X(s,x)-X(s,y)\big],
\end{align*}
\eqref{eq:kolhyp} gives, for every $\theta\in(0,1)$,
\[
\E|G(t)-G(s)|^2\leq4C\min\big(|t-s|^{2\alpha_1},h^{2\alpha_2}\big)\leq4C\,|t-s|^{2\alpha_1\theta}\,h^{2\alpha_2(1-\theta)} .
\]
Let $p\geq1$ and $\eta>0$ with $\eta<\alpha_1\theta-\frac1p$. Since $G$ is Gaussian,
\begin{align*}
\E\int_0^T\!\!\int_0^T\frac{|G(t)-G(s)|^p}{|t-s|^{\eta p+2}}\,dt\,ds
&\leq C_p\,h^{p\alpha_2(1-\theta)}\int_0^T\!\!\int_0^T|t-s|^{p\alpha_1\theta-\eta p-2}\,dt\,ds
=C\,h^{p\alpha_2(1-\theta)},
\end{align*}
the last integral being finite because $p\alpha_1\theta-\eta p-2>-1$. By Step~0 applied to $G$ on $[0,T]$, $\sup_t|G(t)-G(0)|\leq c_{\eta,p}T^\eta B_G^{1/p}$, where $B_G$ is the double integral above. Since $\E|G(0)|^p\leq C_ph^{p\alpha_2}$,
\begin{equation}\label{eq:supG}
\E\Big[\sup_{t\in[0,T]}|X(t,x)-X(t,y)|^p\Big]\leq2^{p-1}\Big(\E|G(0)|^p+c_{\eta,p}^pT^{\eta p}\,\E B_G\Big)\leq C\,|x-y|^{p\alpha_2(1-\theta)} .
\end{equation}
Now let $\gamma_2<\alpha_2$. Choose $\theta$ so small that $\alpha_2(1-\theta)>\gamma_2$, then $p$ so large that $p(\alpha_2(1-\theta)-\gamma_2)>1$, and set
\[
B_x=\int_D\!\!\int_D\frac{\sup_{t\in[0,T]}|X(t,x)-X(t,y)|^p}{|x-y|^{\gamma_2p+2}}\,dx\,dy ,
\]
which is a random variable since the supremum can be taken over rational $t$. By \eqref{eq:supG},
\[
\E B_x\leq C\int_D\!\!\int_D|x-y|^{p\alpha_2(1-\theta)-\gamma_2p-2}\,dx\,dy<\infty .
\]
For each fixed $t$, Step~0 applied to $f_t=X(t,\cdot)$ on $D$, whose double integral is bounded by $B_x$, gives $|X(t,x)-X(t,y)|\leq c_{\gamma_2,p}B_x^{1/p}|x-y|^{\gamma_2}$ for all $x,y\in D$, with a random constant independent of $t$ and with a finite moment of order $p$. Since $p$ can be taken arbitrarily large, this constant has finite moments of all orders.

\emph{Step 3: time increments, uniformly in space.} The argument of Step~2 with the roles of $t$ and $x$ exchanged gives a random constant $K'$ with finite moments of all orders such that $|X(t,x)-X(s,x)|\leq K'|t-s|^{\gamma_1}$ for all $s,t,x$. Taking $K$ as the maximum of the two constants concludes the proof.
\end{proof}

\begin{lemma}[It\^o-Nisio theorem]\label{lem:itonisio}
Let $E$ be a separable Banach space, $(X_k)$ independent symmetric $E$-valued random variables, and $S_n=\sum_{k\leq n}X_k$. If there is an $E$-valued random variable $S$ such that $\langle\phi,S_n\rangle\to\langle\phi,S\rangle$ in probability for every $\phi\in E^*$, then $S_n\to S$ almost surely in $E$.
\end{lemma}

\begin{proof}
Convergence in probability of $\langle\phi,S_n\rangle$ implies the weak convergence of the laws of $\langle\phi,S_n\rangle$ to the law of $\langle\phi,S\rangle$ for every $\phi\in E^*$. This is condition (iv) of \cite[Theorem~2.4]{LT91}, see also \cite{IN68}, which is equivalent to the almost sure convergence of $(S_n)$ to some limit $S'$. Then $\langle\phi,S'\rangle=\langle\phi,S\rangle$ almost surely for every $\phi$ in a countable subset of $E^*$ which separates the points of $E$, hence $S'=S$ almost surely.
\end{proof}

\end{document}